\documentclass[journal]{IEEEtran}
\ifCLASSINFOpdf
   \usepackage[pdftex]{graphicx}
\else
\fi
\usepackage{amsmath}
\usepackage{amssymb}
\usepackage{amsthm}
\usepackage{bm}
\usepackage[compress]{cite}

\usepackage{graphicx}
\usepackage{indentfirst}
\allowdisplaybreaks
\usepackage{subfigure}

\usepackage{lineno}
\usepackage{multirow}
\usepackage{amsfonts}
\usepackage{textcomp}
\usepackage{stfloats}

\usepackage{threeparttable}
\usepackage{hyperref}
\usepackage{float}
\usepackage{indentfirst}
\usepackage{mathrsfs}

\newtheorem{defi}{\textbf{Definition}}
\newtheorem{thom}{\textbf{Theorem}}
\newtheorem{asp}{\textbf{Assumption}}
\newtheorem{rek}{\textbf{Remark}}

\newtheorem{lema}{\textbf{Lemma}}
\newtheorem{algm}{\textbf{Algorithm}}
\newtheorem{pbm}{\textbf{Problem}}
\newtheorem{corol}{\textbf{Corollary}}

\begin{document}
\title{Task Assignment for Multiplayer Reach-Avoid Games in Convex Domains via Analytical Barriers}

\author{Rui~Yan,~\IEEEmembership{Student~Member,~IEEE,}
        Zongying~Shi,~\IEEEmembership{Member,~IEEE,}
        and~Yisheng~Zhong
\thanks{This work was supported by the National Natural Science Foundation of
China under Grants 61374034 and 61210012.}
\thanks{R. Yan, Z. Shi, and Y. Zhong are with the Department of Automation,
Tsinghua University, Beijing 100084, China e-mail: yr15@mails.tsinghua.edu.cn and \{szy, zys-dau\}@mail.tsinghua.edu.cn}
}

\maketitle

\begin{abstract}
This work considers a multiplayer reach-avoid game between two adversarial teams in a general convex domain which consists of a target region and a play region. The evasion team, initially lying in the play region, aims to send as many its team members into the target region as possible, while the pursuit team with its team members initially distributed in both play region and target region, strives to prevent that by capturing the evaders. We aim at investigating a task assignment about the pursuer-evader matching, which can maximize the number of the evaders who can be captured before reaching the target region safely when both teams play optimally. To address this, two winning regions for a group of pursuers to intercept an evader are determined by constructing an analytical barrier which divides these two parts. Then, a task assignment to guarantee the most evaders intercepted is provided by solving a simplified 0-1 integer programming instead of a non-deterministic polynomial problem, easing the computation burden dramatically. It is worth noting that except the task assignment, the whole analysis is analytical. Finally, simulation results are also presented.
\end{abstract}

\begin{IEEEkeywords}
Reach-avoid games; Multi-agent systems; Barriers; Pursuit-evasion games; Optimal control; Differential games 
\end{IEEEkeywords}

\section{Introduction}
\IEEEPARstart{T}{his} paper studies a multiplayer reach-avoid (RA) game between two adversarial teams of cooperative players playing in a convex bounded planar domain, which is partitioned into a target region and a play region by a straight line. Starting from the play region, the evasion team aims to send as many its team members, called evaders, into the target region as possible. Conversely, the pursuit team initially lying in the target region and play region, strives to prevent the evasion team from doing so by attempting to capture the evaders. Actually, from another side, this game can also be viewed as an evasion team tries to escape from a bounded region through an exit which is represented by a straight line, while avoiding adversaries and moving obstacles formulated as a pursuit team. Such a differential game is a powerful theoretical tool for analyzing realistic situations in robotics, aircraft control, security, reachability analysis and other domains \cite{Ba1999Dynamic,Chung2011,shaferman2017cooperative,Dong7862774,DING20132665}. For example, in collision avoidance and path planning, how a group of vehicles can get into some target set or escape from a bounded region through an exit, while avoiding dangerous situations, such as collisions with static or moving obstacles \cite{Zhou6426643,Panagou6767157,Pan2012Pursuit}. In region pursuit games, multiple pursuers are used to intercept multiple adversarial intruders \cite{ChenMo2016Multiplayer,HAYOUN2017On,Raslan2016A}. In safety verification, an agent often needs to judge whether it can guarantee its arrival into a safe region throughout plenty of dynamic dangers, such as disturbances and adversaries \cite{Huang2015Aut}. In the references \cite{Scott6580287} and \cite{Garcia8340791}, cooperative behaviors within pursuit-evasion games are analyzed in order to help or rescue teammates in the presence of adversarial players.

However, finding cooperative strategies and task assignment among multiple players, especially when two teams both consist of multiple members, can be challenging \cite{Chung8424838}, as computing solutions over the joint state space of multiple players can greatly increase computational complexity, beyond the scope of traditional dynamic programming \cite{Ba1999Dynamic}. In \cite{Gerkey04045564}, a formal analysis and taxonomy of task allocation in multi-robot systems is presented. More recently, a distributed version of the Hungarian method \cite{Chopra7932518} and consensus-based decentralized auction and bundle algorithms \cite{Choi5072249} are proposed to solve the multirobot assignment problem. The authors in \cite{Zhou5735231} study the problem of multirobot active target tracking by processing relative observations. Due to the conflicting and asymmetric goals between two teams, complex cooperations and non-intuitive strategies within each team may exist. Although some techniques have been used to analyze RA games and demonstrated to work well in some conditions, they are mostly numerical algorithms and suffer from some weaknesses, such as highly computational complexity, conservation, strong assumption and application limitations \cite{Xue2017Reach,Chen8267187,Vamvoudakis2011Multi,ZHOU201828}. 

Although the problem considered in this paper is different from the classical pursuit-evasion games that have been thoroughly studied, we borrow several existing notions and modify their definitions slightly to address our current scenario. For RA games, as Isaacs' book \cite{Is1967Diff} shows, the core point is to construct the barrier, which is the boundary of the RA set, splitting the entire state space into two disjoint parts: Pursuit Winning Region (PWR) and Evasion Winning Region (EWR). The PWR is the region of initial conditions, from which the pursuit team can ensure the capture before the evader enters the target region. The EWR, complementary to the PWR, is the region of initial conditions, from which the evader can succeed to reach the target region regardless of the pursuit team's strategies. The surface that separates the PWR from the EWR is called barrier. 

In principle, the Hamilton-Jacob-Isaacs (HJI) approach is an ideal tool for solving general RA games when the game is low-dimensional, such as autonomous river navigation for underactuated vehicles \cite{Weekly6838996} and safety specifications in hybrid systems \cite{Tomlin2000A}. By defining a value function merging the payoff function and discriminator function with minmax operation, this approach involves solving a HJI partial differential equation (PDE) in the joint state space of the players and locating the barrier by finding the zero sublevel set of this value function. The players' optimal strategies can be extracted from the gradient of the value function. Generally, there are two approaches to solve the HJI PDEs: the method of characteristics \cite{Garcia8340791} and numerical approximation of the value function on a grid of the continuous state space \cite{Mitchell2005A,Margellos2011Ham,Fis2015Rea}.

However, in practical applications, two approaches both face computational challenges. Non-unique terminal conditions in RA game setups, capture or entry into the target set, make it difficult to generate strategies by characteristic solutions which require backward integration from terminal manifold, as different backward trajectories may produce complicated singular surfaces for which there exist no systematic analysis methods \cite{lewin2012differential}. On the other hand, a number of numerical tools for solving HJI PDEs on grids have been provided to solve practical problems \cite{Margellos2011Ham}. Unfortunately, the curse of dimensionality makes these approaches computationally intractable in our multiplayer RA games, as the grid required for approximating the value function scales exponentially with the number of players.

For certain games and game setups, geometric method shows an incredible power in providing strategies for the players \cite{BOPARDIKAR20091771,Yan2017Escape,lopez2016optimal}. For example, Voronoi diagrams, dividing a plane into regions of points that are closest to a predetermined set of seed points, are widely used for generating strategies in pursuit-evasion games, usually when each player possesses the same speed. Especially in group pursuit of a single evader or multiple evaders, Voronoi-based approaches can provide very constructive cooperative strategies, such as minimizing the area of the generalized Voronoi partition of the evader \cite{Pierson2017Int,Zhou2016Cooperative} or pursuing the evader in a relay way \cite{Bakolas2012Relay}. As for unequal speed scenarios, the Apollonius circle, first introduced by Isaacs, is a useful tool for analyzing the capture of a high-speed evader by using multiple pursuers \cite{Ramana2016Pur,Yan2017Rea}. More realistically, when pursuit-evasion games are played in the presence of obstacles that inhibit the motions of the players, Euclidean shortest path method is employed to construct the dominance region in \cite{Oyler2016Pursuit}, and visibility-based target tracking games are addressed in \cite{Zou8588382} and \cite{Yu6112246}. The work by Katsev \emph{et al.} \cite{Katsev5685659} introduces a simple wall-following robot to map and solve pursuit-evasion strategies in an unknown polygonal environment. The authors \cite{Noori913498291} revisit the lion and man problem by introducing line-of-sight visibility. In \cite{Bhadauria4912452894}, the number of pursuers which can guarantee to capture an equal speed evader in polygonal environment with obstacles, is investigated. For RA games, a number of straight lines called paths of defense, separating the target set from evaders, are created successively to compute an approximate 2D slice of the reach-avoid set, which eases the computational burden sharply \cite{ChenMo2016Multiplayer}.

The construction of barrier, the most central and important part in RA games, has attracted a lot of attention in pursuit-evasion games and until now, achieved remarkable results \cite{ruiz2016differential,Sun2015Pur,Zou2016Vis,Yan2017D}. For example, in \cite{Ruiz6525413} and \cite{MACIAS2018271}, the authors compute the barrier for a pursuit-evasion game between an omnidirectional evader and a differential drive robot. For the problem of tracking an evader in an environment containing a corner, the method of explicit policy is used to investigate the escape set and the track set \cite{Sourabh11415885}. In the reference \cite{Barron2018}, the boundary of reach-avoid set is studied for general reach-avoid differential games. By adopting the method of characteristics or numerical approximation on grid, the barrier for two players or at most three players is constructed analytically or numerically. However, computing the barrier directly for more than three players in pursuit-evasion games is missing, resulting from the intrinsic high dimensionality of the joint state space.

Compared with the traditional pursuit-evasion games, RA games are more complicated and have more practical significance, as the evaders aim to not only avoid the capture, but also strive to reach a target set. To our best knowledge, the current work for RA games mainly focuses on the construction of barrier in two-player scenarios \cite{ChenMo2016Multiplayer,Margellos2011Ham,Fis2015Rea}. The methods in \cite{ChenMo2016Multiplayer,Margellos2011Ham,Fis2015Rea} are numerical and cannot directly obtain the barrier for multiple players due to the high computation burden. The work \cite{Yan2017Rea} is most similar to this work, but it only considers two pursuers and one evader case without involving the complicated cooperations among the pursuit team occurring in this work, and it only focuses on a square game domain which is quite limited. Moreover, our current work also allows the pursuers to start the game in the target region, which is not considered in \cite{Yan2017Rea}.

In this work, an analytical study on the number of the evaders which the pursuit team would be able to prevent from reaching the target region, is presented by generating a task assignment, namely, pursuer-evader matching pairs. Actually, the above analysis refers to a game of kind \cite{Is1967Diff}. To this end, all pursuers are first classified into all possible coalitions. Then, an evasion region method is proposed to construct the barrier analytically for each pursuit coalition versus one evader in convex domains, which is the first time in the existing literature to construct the barrier directly for multiplayer games with more than three players involved. More importantly, the constructed barrier is analytical and overcomes the curse of dimensionality. Finally, with rich prior information in hand about which evaders can be intercepted by a specified pursuit coalition, a maximum pursuer-evader matching is given such that the most evaders are intercepted.

The original contributions of this paper are as follows. First, the analytical barrier is constructed for one pursuer and two pursuers versus one evader in convex domains. Second, the analytical barrier for multiple pursuers versus one evader is given, involving two kinds of initial deployments shown in Assumptions \ref{constrainedinitial} and \ref{relaxedinitial}, which can determine the capturable and uncapturable regions for these pursuers. Third, since all possible cooperations among the pursuit team are considered, the upper bound on the number of the evaders which the pursuit team can guarantee to intercept, is given by solving a 0-1 integer programming instead of a non-deterministic polynomial problem, greatly easing the computation burden. Fourth, except the task assignment, the whole analysis is analytical, allowing for real-time updates. These contributions provide a complete solution from the perspective of the task assignment to multiplayer RA games in any convex domains consisting of a target region and a play region separated by a straight line.

The rest of this paper is organized as follows. In Section \ref{problmstate}, the problem statement is given. In Section \ref{preliminaries}, some important preliminaries are presented. In Section \ref{pursuitcoalition}, the barrier and winning regions for one pursuit coalition versus one evader, are found. In Section \ref{task}, the task assignment for the pursuit team to capture the most evaders, is designed. In Section \ref{simulation}, simulation results are presented. Finally, Section \ref{conclusion} concludes the paper and the future work is discussed.

\begin{figure}\centering
\graphicspath{{Figures/}}
\includegraphics[width=58mm,height=43mm]{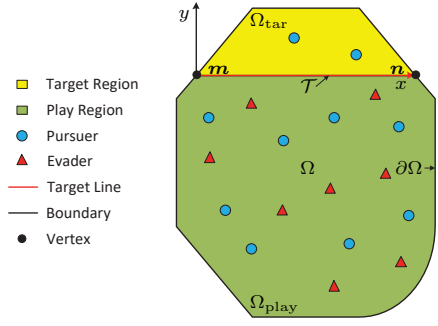}
\put(-72,5){$\scriptstyle{\Omega_{\rm play}}$}
\put(-72,112){$\scriptstyle{\Omega_{\rm tar}}$}
\put(-99,115){$\scriptstyle{y}$}
\put(-17,88){$\scriptstyle{x}$}
\put(-88,94){$\scriptstyle{\bm{m}}$}
\put(-19,94){$\scriptstyle{\bm{n}}$}
\put(-53,87){$\scriptstyle{\mathcal{T}}$}
\put(-17,56){$\scriptstyle{\partial\Omega}$}
\put(-53,56){$\scriptstyle{\Omega}$}
\caption{Multiplayer reach-avoid games in convex domains, where the pursuit team with multiple pursuers (blue circles) wants to capture the most evaders (red triangles) before these evaders enter the target region. Our goal is to find a task assignment for the pursuit team to guarantee the most evaders captured, involving pursuer-evader matching pairs.\label{fig:1}}
\end{figure}

\section{Problem Statement}\label{problmstate}
\subsection{Multiplayer Reach-Avoid Games}
Consider $N_p+N_e$ players partitioned into two teams, a pursuit team of $N_p$ pursuers, $\{P_i\}_{i=1}^{N_p}=\{{P_1},...,{P_{N_p}}\}$, and  an evasion team of $N_e$ evaders, $\{E_j\}_{j=1}^{N_e}=\{{E_1},...,{E_{N_e}}\}$, whose states are constrained in a bounded, convex domain $\Omega\subset \mathbb{R}^2$ with boundary $\partial\Omega$. Each player is assumed to be a mass point. As Fig. \ref{fig:1} shows, a straight line $\mathcal{T}\subset\Omega$, called target line, divides $\Omega$ into two disjoint parts $\Omega_{\rm tar}$ and $\Omega_{\rm play}$. The compact set $\Omega_{\rm tar}$, called target region, represents the region which the evasion team strives to enter while the pursuit team tries to protect. The play region $\Omega_{\rm play}=\Omega\setminus\Omega_{\rm tar}$ corresponds to the region in which two teams play the game. Also note that $\mathcal{T}\subset\Omega_{\rm tar}$. Let $\mathbf{x}_{P_i}(t)=({x}_{P_i}(t),y_{P_i}(t))\in\mathbb{R}^2$ and $\mathbf{x}_{E_j}(t)=({x}_{E_j}(t),y_{E_j}(t))\in\mathbb{R}^2$ be the positions of $P_i$ and $E_j$ at time $t$, respectively. The dynamics of the players are described by the following decoupled system for $t\ge0$:
\begin{equation}\begin{aligned}
 \dot{\mathbf{x}}_{P_i}(t)&=v_{P_i}\mathbf{u}_{P_i}(t),& \hspace{0.5mm} \mathbf{x}_{P_i}(0)&=\mathbf{x}^0_{P_i},& i&=1,...,N_p
 \\\dot{\mathbf{x}}_{E_j}(t)&=v_{E_j}\mathbf{u}_{E_j}(t),&\mathbf{x}_{E_j}(0)&=\mathbf{x}^0_{E_j},& j&=1,...,N_e
\end{aligned} \end{equation}
where $\mathbf{x}^0_{P_i}=({x}_{P_i}^0,{y}_{P_i}^0)\in\mathbb{R}^2$ and $\mathbf{x}^0_{E_j}=({x}_{E_j}^0,{y}_{E_j}^0)\in\mathbb{R}^2$ are the initial positions of $P_i$ and $E_j$, respectively. The maximum speeds of $P_i$ and $E_j$ are denoted by $v_{P_i}$ and $v_{E_j}$ respectively. The control inputs at time $t$ for $P_i$ and $E_j$ are $\mathbf{u}_{P_i}(t)$ and $\mathbf{u}_{E_j}(t)$ respectively, and they satisfy the constraint $\mathbf{u}_{P_i}(t),\mathbf{u}_{E_j}(t)\in\mathcal{U}=\{\mathbf{u}\in\mathbb{R}^2|\|\mathbf{u}\|_2=1\}$, where $\|\cdot\|_2$ stands for the Euclidean norm in $\mathbb{R}^2$. Unless needed for clarity, to simplify notations, $t$ will be omitted hereinafter. 

The goal of the evasion team is to send as many evaders as possible into $\Omega_{\rm tar}$ without being captured, while the pursuit team strives to prevent the evasion team from that by capturing the evaders. Naturally, once an evader enters $\Omega_{\rm tar}$, the pursuer cannot capture it hereinafter. Assume that $E_j$ is captured by $P_i$ in $\Omega_{\rm play}$ if $E_j$'s position coincides with $P_i$'s position, that is, the point-capture is considered. Our paper aims to investigate an optimal interception matching scheme for the pursuit team such that the most evaders can be intercepted, and also provide strategy instructions for the evasion team.

\begin{figure}\flushright
\graphicspath{{Figures/}}
\includegraphics[width=61mm,height=6mm]{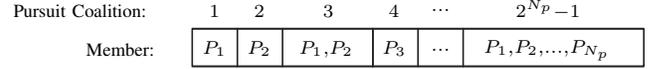}
\put(-240,20){\scriptsize{Pursuit Coalition:}}
\put(-212.5,5){\scriptsize{Member:}}
\put(-168.5,6){$\scriptstyle{{P_1}}$}
\put(-166,20){$\scriptstyle{1}$}
\put(-149,20){$\scriptstyle{2}$}
\put(-123,20){$\scriptstyle{3}$}
\put(-98,20){$\scriptstyle{4}$}
\put(-50,20){$\scriptstyle{2^{N_p}-1}$}
\put(-151.5,6){$\scriptstyle{{P_2}}$}
\put(-132,6){$\scriptstyle{{P_1},{P_2}}$}
\put(-100.5,6){$\scriptstyle{{P_3}}$}
\put(-82,7){$\scriptstyle{...}$}
\put(-82,22){$\scriptstyle{...}$}
\put(-62,6.3){$\scriptstyle{{P_1},{P_2},...,{P_{N_p}}}$}
\caption{All possible pursuit coalitions in a pursuit team with $N_p$ pursuers coded in a binary way. For example, the binary code of 3 is $11$, and thus the pursuit coalition 3's members are $P_1$ and $P_2$, namely, $m_1=1$ and $m_2=2$.\label{fig:4}}
\end{figure}

In view of the pursuit team's goal, all possible cooperations among the pursuers must be considered, which is not involved in the existing literature on multiplayer RA games. The maximum matching employed by \cite{ChenMo2016Multiplayer} can only provide a suboptimal strategy for the pursuit team, as cooperation is only introduced at the matching step. Since the pursuit team has $N_p$ pursuers, one can select one-pursuer coalitions, two-pursuer coalitions...until $N_p$-pursuer coalitions, as Fig. \ref{fig:4} shows. Therefore, there are $2^{N_p}-1$ alternatives for pursuit coalitions among $N_p$ pursuers. These coalitions are labeled from 1 to $2^{N_p}-1$ successively, so that they can be coded in a binary way as follows.

\begin{defi}[Binary Coalitions for Pursuit Team]\label{binary coalitions}\rm
$\forall k=1,...,2^{N_p}-1$, $P_i$ belongs to the pursuit coalition $k$ if the $i$-th bit from low order side in binary representation of $k$ is 1. Then, denote the index set of the pursuers in pursuit coalition $k$ by $\mathcal{I}_k=\{m_j|1\leq m_j\leq N_p,j=1,2,...,n_k\}$ satisfying $m_{i}<m_{j}$ for $i<j$, where $n_k$ is the number of the pursuers in pursuit coalition $k$.
\end{defi}

\begin{defi}[Pursuit Subcoalition]\label{Pursuit Subcoalition}\rm
Consider two pursuit coalitions $k_1$ and $k_2$. If every pursuer in $k_1$ occurs in $k_2$, then $k_1$ is called a pursuit subcoalition of $k_2$.
\end{defi}

Obviously, there are plenty of ways to code these pursuit coalitions, but this binary way is very convenient to determine every coalition's members by only recording its number, as Fig. \ref{fig:4} shows. For example, for the pursuit coalition $k=5$, since the binary representation of $5$ is $101$, thus this pursuit coalition's members are $P_1$ and $P_3$, namely, $m_1=1$ and $m_2=3$. The adoption of pursuit subcoalition will tremendously simplify our problems as discussed below. 

\subsection{Information Structure and Assumptions}\label{assumtionsection}
As is the usual convention in the differential game theory, the equilibrium outcomes crucially depend on the information structure employed by each player. Classically, the state feedback information structure allows each player to choose its current input, $\mathbf{u}_{P_i}$ or $\mathbf{u}_{E_j}$, based on the current value of the information set $\{\mathbf{x}_{P_1},\cdots,\mathbf{x}_{P_{N_p}},\mathbf{x}_{E_1},\cdots,\mathbf{x}_{E_{N_e}}\}$. This paper focuses on a non-anticipative information structure, as commonly adopted in the differential game literature (see for example \cite{Mitchell2005A}, \cite{Elliott1972The}). Under this information structure, the pursuit team is allowed to make decisions about its current input with all the information of state feedback, plus the evasion team's current input. While the evasion team is at a slight disadvantage under this information structure, at a minimum he has access to sufficient information to use state feedback, because the pursuit team must declare his strategy before the evasion team chooses a specific input and thus the evasion team can determine the response of the pursuit team to any input signal. Thus, the multiplayer reach-avoid games formulated here are an instantiation of the Stackelberg game \cite{Ba1999Dynamic}.

As Fig. \ref{fig:1} shows, let $\bm{m}$ and $\bm{n}$ denote the endpoints of $\mathcal{T}$, and assume $\|\bm{m}-\bm{n}\|_2=l$. Fix the origin at $\bm{m}$, and build a Cartesian coordinate system with $x$-axis along the straight line through $\bm{m}$ and $\bm{n}$, and $y$-axis perpendicular to $x$-axis and pointing to $\Omega_{\rm tar}$. For unity, we assume that $\mathbb{R}^n=\mathbb{R}^{1\times n}$, where $n$ is a positive integer.

Next, it is assumed that the following conditions are satisfied by the initial configurations of the players, where Assumptions \ref{constrainedinitial} and \ref{relaxedinitial} will be separately considered.

\begin{asp}[Isolate Initial Deployment]\label{isolateinitial}\rm
The initial positions of the players satisfy the three conditions:\\
1) $\|\mathbf{x}^0_{P_i}-\mathbf{x}^0_{P_j}\|_2>0$ for all $i,j=1,...,N_p,i\neq j$;\\
2) $\|\mathbf{x}^0_{E_i}-\mathbf{x}^0_{E_j}\|_2>0$ for all $i,j=1,...,N_e,i\neq j$;\\
3) $\|\mathbf{x}^0_{P_i}-\mathbf{x}^0_{E_j}\|_2>0$ for all $i=1,...,N_p,j=1,...,N_e$.
\end{asp}

It can be seen that \autoref{isolateinitial} guarantees that all players start to play the game from different initial positions and every evader is not captured by the pursuers initially. 

\begin{asp}[Constrained Initial Deployment]\label{constrainedinitial}\rm
Suppose that $\mathbf{x}^0_{P_i}\in\Omega_{\rm play}\cup\mathcal{T}$ for all $i=1,...,N_p$ and $\mathbf{x}^0_{E_j}\in\Omega_{\rm play}$ for all $j=1,...,N_e$.
\end{asp}

\begin{asp}[Relaxed Initial Deployment]\label{relaxedinitial}\rm
Suppose that $\mathbf{x}^0_{P_i}\in\Omega$ for all $i=1,...,N_p$ and $\mathbf{x}^0_{E_j}\in\Omega_{\rm play}$ for all $j=1,...,N_e$.
\end{asp}

In Assumptions 2 and 3, restricting the evaders' initial positions into $\Omega_{\rm play}$ is a reasonable assumption for our RA games, as the evader wins once it enters $\Omega_{\rm tar}$. As for the initial positions of the pursuers, Assumption 2 is employed out of consideration for developing an extensible basic approach. Then, the case with more practical initial configuration described in Assumption 3 is investigated by building a bridge to the proposed basic approach. Moreover, an initial deployment is admissible if Assumptions 1 and 2 or 1 and 3 hold. 

Generally, in multiplayer RA games, the pursuers are homogeneous and the same for the evaders, such as confrontation between two species, and collision avoidance in the environment with similar dynamic obstacles. Thus, our discussion assumes that the pursuers have the same maximum speed $v_P$, and the evaders have the same maximum speed $v_E$. Define $\alpha=v_E/v_P$ to be the speed ratio, and this paper focuses on the faster pursuers.


\begin{asp}[Speed Ratio]\label{speedratio}\rm
Suppose that $v_{P_i}=v_P>0,i=1,...,N_p$, and $v_{E_j}=v_E>0,j=1,...,N_e$. Assume that $\alpha=v_E/v_P$ satisfies $0<\alpha<1$.
\end{asp}


\begin{figure}\centering
\graphicspath{{Figures/}}
\includegraphics[width=78mm,height=24mm]{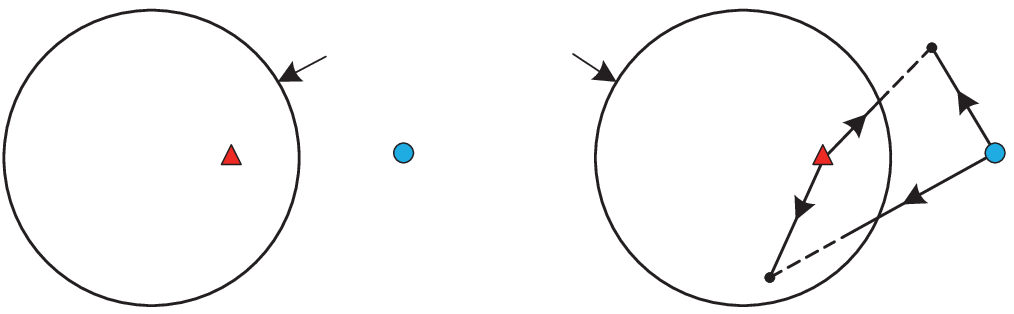}
\put(-54,39){$\scriptstyle{E_j}$}
\put(-105,57){$\scriptstyle{{\mathcal{A}}}$}
\put(-138,26){$\scriptstyle{P_i}$}
\put(-9,26){$\scriptstyle{P_i}$}
\put(-176,26){$\scriptstyle{E_j}$}
\put(-150,56){$\scriptstyle{{\mathcal{A}}}$}
\put(-194,32){$\scriptstyle{\mathcal{R}_e}$}
\put(-70,32){$\scriptstyle{\mathcal{R}_e}$}
\put(-55,-8){$\scriptstyle{(b)}$}
\put(-180,-8){$\scriptstyle{(a)}$}
\caption{The evasion region (ER) and the boundary of ER (BER). $(a)$ The ER and BER determined by $P_i$ and $E_j$ are $\mathcal{R}_e$ and ${\mathcal{A}}$ respectively, where ${\mathcal{A}}$ is a circle and $\mathcal{R}_e$ is the interior of this circle. $(b)$ The circle ${\mathcal{A}}$, also called Apollonius circle, divides $\mathbb{R}^2$ into two parts: The interior of the circle is $E_j$'s dominance region, i.e. the ER $\mathcal{R}_e$: $E_j$ can reach any point inside the circle before $P_i$; the exterior of the circle is $P_i$'s dominance region: $P_i$ can reach any point outside the circle before $E_j$.\label{figT:1}}
\end{figure}

\section{Preliminaries}\label{preliminaries}
\subsection{Computation of the ER and BER}
For $\bm{x},\bm{y}\in\mathbb{R}^2$, define two maps $r(\bm{x},\bm{y})=\frac{\alpha\|\bm{x}-\bm{y}\|_2}{1-\alpha^2}$ and $\eta(\bm{x},\bm{y})=\frac{\bm{x}-\alpha^2\bm{y}}{1-\alpha^2}$ \cite{Ramana2016Pur}, whose geometric meanings will be stated below.

Let the set of points in $\mathbb{R}^2$ that one evader can reach before one pursuer, regardless of the pursuer' best effort, be called evasion region (ER), and the surface which bounds ER is called the boundary of ER (BER). 

Denote the ER and BER determined by $P_i$ and $E_j$ by $\mathcal{R}_e(\mathbf{x}_{E_j}^0,\mathbf{x}_{P_i}^0)$ and ${\mathcal{A}}(\mathbf{x}_{E_j}^0,\mathbf{x}_{P_i}^0)$ respectively, which can be mathematically formulated as follows:
\begin{equation}\begin{aligned}\label{ARBAR1v1}
\mathcal{R}_e&=\big\{\mathbf{z}\in\mathbb{R}^2|\|\mathbf{z}-\mathbf{x}_{E_j}^0\|_2<\alpha\|\mathbf{z}-\mathbf{x}_{P_i}^0\|_2\big\}\\
{\mathcal{A}}&=\big\{\mathbf{z}\in\mathbb{R}^2|\|\mathbf{z}-\mathbf{x}_{E_j}^0\|_2=\alpha\|\mathbf{z}-\mathbf{x}_{P_i}^0\|_2\big\}.
\end{aligned}\end{equation}
Also note that\begin{equation}\begin{aligned}\label{apolloniuscircle22}
\|\mathbf{z}-\mathbf{x}_{E_j}^0\|^2_2<\alpha^2\|\mathbf{z}-\mathbf{x}_{P_i}^0\|^2_2\Rightarrow\Big\|\mathbf{z}-\frac{\mathbf{x}_{E_j}^0-\alpha^2\mathbf{x}_{P_i}^0}{1-\alpha^2}\Big\|_2^2\qquad&\\
<\frac{\alpha^2\|\mathbf{x}_{E_j}^0-\mathbf{x}_{P_i}^0\|^2_2}{(1-\alpha^2)^2}\Rightarrow\|\mathbf{z}-\eta(\mathbf{x}_{E_j}^0,\mathbf{x}_{P_i}^0)\|_2^2<r^2(\mathbf{x}_{E_j}^0,\mathbf{x}_{P_i}^0)&
\end{aligned}\end{equation}
implying that (\ref{ARBAR1v1}) can be equivalently rewritten as
\begin{equation}\begin{aligned}\label{apolloniuscircle33}
\mathcal{R}_e&=\big\{\mathbf{z}\in\mathbb{R}^2|\|\mathbf{z}-\eta(\mathbf{x}_{E_j}^0,\mathbf{x}_{P_i}^0)\|_2<r(\mathbf{x}_{E_j}^0,\mathbf{x}_{P_i}^0)\big\}\\
{\mathcal{A}}&=\big\{\mathbf{z}\in\mathbb{R}^2|\|\mathbf{z}-\eta(\mathbf{x}_{E_j}^0,\mathbf{x}_{P_i}^0)\|_2=r(\mathbf{x}_{E_j}^0,\mathbf{x}_{P_i}^0)\big\}.
\end{aligned}\end{equation}
Thus, it can be seen that ${\mathcal{A}}$ is a circle of radius $r(\mathbf{x}_{E_j}^0,\mathbf{x}_{P_i}^0)$ centered at $\eta(\mathbf{x}_{E_j}^0,\mathbf{x}_{P_i}^0)$, and $\mathcal{R}_e$ is the interior of ${\mathcal{A}}$, shown in Fig. \ref{figT:1}(a). This circle ${\mathcal{A}}$, also called Apollonius circle \cite{Is1967Diff} divides $\mathbb{R}^2$ into two parts: The interior of the circle is $E_j$'s dominance region, i.e., the ER $\mathcal{R}_e$: $E_j$ can reach any point inside the circle before $P_i$; the exterior of the circle is $P_i$'s dominance region: $P_i$ can reach any point outside the circle before $E_j$, as Fig. \ref{figT:1}(b) illustrates. 


Unless needed for clarity, to simplify notations, drop the initial positions occurring in the expressions of $\mathcal{R}_e$ and ${\mathcal{A}}$. 
 
\subsection{Key Function}
Next, present an important lemma about a key function.
\begin{lema}[Monotony of Function]\label{monotonic function1}\rm
Given $P_i$ and $E_j$'s initial positions $\mathbf{x}_{P_i}^0$ and $\mathbf{x}_{E_j}^0$ satisfying $y_{E_j}^0<0$ and $x_{P_i}^0\leq x_{E_j}^0$, if ${\mathcal{A}}\cap\{\mathbf{z}\in\mathbb{R}^2|y=0\}=\big\{(c_1,0),(c_2,0)\big\}$ with $c_1<c_2$, the function
\begin{equation}\label{monotonicfunction1}
G_1(x_p)=\|\bm{p}-\mathbf{x}_{P_i}^0\|_2-\frac{\|\bm{p}-\mathbf{x}_{E_j}^0\|_2}{\alpha},\bm{p}=(x_p,0)
\end{equation}
is strictly monotonic increasing when $x_p\in[c_1,x_p^*]$, and strictly monotonic decreasing when  $x_p\in[x_p^*,c_2]$, where $x_p^*$ is the unique solution of the quartic equation
\begin{equation}\label{quarticlemma}
\frac{x_p^*-x_{P_i}^0}{\|\bm{p}^*-\mathbf{x}_{P_i}^0\|_2}=\frac{x_p^*-x_{E_j}^0}{\alpha\|\bm{p}^*-\mathbf{x}_{E_j}^0\|_2},\bm{p}^*=(x_p^*,0)
\end{equation}
in the interval $[c_1,c_2]$.

\emph{Proof:} See Fig. \ref{figT:2}(a). Take $\bm{p}=(x_p,0)$. It can be seen that if $x_p\in[c_1,c_2]$, $G_1(x_p)$ in (\ref{monotonicfunction1}) is actually the distance (depicted in dashed line) between $P_i$ and $E_j$ exactly when $E_j$ arrives at $\bm{p}$, if $P_i$ and $E_j$ both move directly towards $\bm{p}$.

Note that ${\mathcal{A}}$ is a circle. Thus, based on (\ref{ARBAR1v1}), ${\mathcal{A}}\cap\{\mathbf{z}\in\mathbb{R}^2|y=0\}=\big\{(c_1,0),(c_2,0)\big\}$ with $c_1<c_2$ implies that $G_1(c_1)=G_1(c_2)=0$, $G_1(x_p)>0$ for $x_p\in(c_1,c_2)$, and $G_1(x_p)<0$ for $x_p\in(-\infty,c_1)\cup(c_2,+\infty)$. Thus, the maximum point $x_p^*$ of $G_1(x_p)$ lies in $[c_1,c_2]$ and satisfies $G_1'(x_p^*)=0$. If we can verify $G_1'(x_p^*)=0$ admits a unique solution $x_p^*$ in the interval $[c_1,c_2]$, then the lemma is straightforward. First, the existence of $x_p^*$ in the interval $[c_1,c_2]$ is obvious by noting that $G_1(x_p)$ is continuous in this interval. Next, prove the uniqueness.

\begin{figure}\centering
\graphicspath{{Figures/}}
\includegraphics[width=75mm,height=25.6mm]{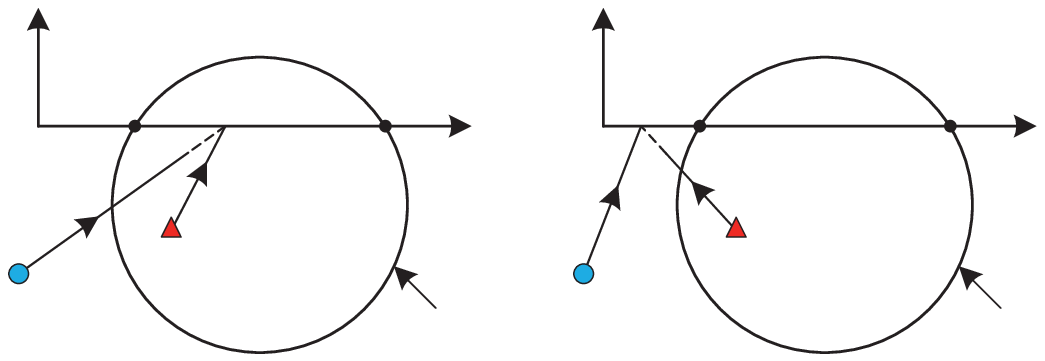}
\put(-67,18){$\scriptstyle{E_j}$}
\put(-98,9){$\scriptstyle{P_i}$}
\put(-10,42){$\scriptstyle{x}$}
\put(-96,62){$\scriptstyle{y}$}
\put(-127,5){$\scriptstyle{{\mathcal{A}}}$}
\put(-214,9){$\scriptstyle{P_i}$}
\put(-183,18){$\scriptstyle{E_j}$}
\put(-169,50){$\scriptstyle{\bm{p}}$}
\put(-83.5,50){$\scriptstyle{\bm{p}}$}
\put(-204,52){$\scriptstyle{(c_1,0)}$}
\put(-70,40){$\scriptstyle{(c_1,0)}$}
\put(-134.5,52){$\scriptstyle{(c_2,0)}$}
\put(-18.5,52){$\scriptstyle{(c_2,0)}$}
\put(-126,42){$\scriptstyle{x}$}
\put(-212,63){$\scriptstyle{y}$}
\put(-11,5){$\scriptstyle{{\mathcal{A}}}$}
\put(-55,-8){$\scriptstyle{(b)}$}
\put(-170,-8){$\scriptstyle{(a)}$}
\caption{The monotony of the function $G_1(x_p)$, where $\bm{p}=(x_p,0)$ and $G_1(c_1)=G_1(c_2)=0$. $(a)$ If $x_p\in(c_1,c_2)$, $G_1(x_p)>0$, that is, $E_j$ can reach $\bm{p}$ before $P_i$. $(b)$ If $x_p\in(-\infty,c_1)\cup(c_2,+\infty)$, $G_1(x_p)<0$, that is, $P_i$ can reach $\bm{p}$ before $E_j$. It is proved that there exists a unique maximum point of $G_1(x_p)$ for $x_p\in[c_1,c_2]$.\label{figT:2}}
\end{figure}

Take $\bm{p}^*=(x_p^*,0)$. By taking the derivative for (\ref{monotonicfunction1}) with respect to $x_p$, $G_1'(x_p^*)=0$ means that (\ref{quarticlemma}) holds. There are two cases depending on whether $x_{P_i}^0<x_{E_j}^0$ or $x_{P_i}^0=x_{E_j}^0$, which will be separately discussed below.

Consider $x_{P_i}^0<x_{E_j}^0$ first. Since $x_p^*\in[c_1,c_2]$, as Fig. \ref{figT:2}(a) shows, we have $G_1(x_p^*)\ge0$, implying that
\begin{equation}\label{lemmaG1de344}
\|\bm{p}^*-\mathbf{x}_{P_i}^0\|_2\ge\frac{\|\bm{p}^*-\mathbf{x}_{E_j}^0\|_2}{\alpha}>\alpha\|\bm{p}^*-\mathbf{x}_{E_j}^0\|_2.
\end{equation}
Then, combining (\ref{quarticlemma}) and (\ref{lemmaG1de344}) leads to
\begin{equation}\label{lemmaG1de3554}
|x_p^*-x_{P_i}^0|>|x_p^*-x_{E_j}^0|.
\end{equation}
Also note that (\ref{quarticlemma}) guarantees that $x_p^*-x_{P_i}^0$ and $x_p^*-x_{E_j}^0$ have the same plus or minus sign. Thus, by also noting that $x_{P_i}^0<x_{E_j}^0$, then (\ref{lemmaG1de3554}) means that $x_p^*>x_{E_j}^0>x_{P_i}^0$, as Fig. \ref{figT:2}(a) illustrates. Define
\begin{equation}\label{lemmaG1de3}
G_2(x_p^*)=\frac{\|\bm{p}^*-\mathbf{x}_{P_i}^0\|_2^2(x_p^*-x_{E_j}^0)^2}{\|\bm{p}^*-\mathbf{x}_{E_j}^0\|_2^2(x_p^*-x_{P_i}^0)^2}
\end{equation}
and thus (\ref{quarticlemma}) can be rewritten as $G_2(x_p^*)=\alpha^2$. Since $\alpha^2<1$, then $G_2(x_p^*)<1$, implying that
\begin{equation}\begin{aligned}\label{lemmaG1de4}
&\|\bm{p}^*-\mathbf{x}_{E_j}^0\|_2^2(x_p^*-x_{P_i}^0)^2>\|\bm{p}^*-\mathbf{x}_{P_i}^0\|_2^2(x_p^*-x_{E_j}^0)^2\\
&\Rightarrow(y_{E_j}^0)^2(x_p^*-x_{P_i}^0)^2-(y_{P_i}^0)^2(x_p^*-x_{E_j}^0)^2>0.
\end{aligned}\end{equation}
By computing the derivative of $G_2(x_p^*)$ with respect to $x_p^*$, it can be verified that $G_2(x_p^*)$ is strictly monotonic for $x_p^*$ satisfying (\ref{lemmaG1de4}) and $x_p^*>x_{E_j}^0>x_{P_i}^0$. Thus, $G_2(x_p^*)=\alpha^2$, i.e., $G'_1(x_p^*)=0$, admits a unique solution $x_p^*$ in the interval $[c_1,c_2]$.

If $x_{P_i}^0=x_{E_j}^0$, (\ref{quarticlemma}) implies that $x_p^*=x_{E_j}^0=x_{P_i}^0$ by noting (\ref{lemmaG1de344}). Thus, $G'_1(x_p^*)=0$ still admits a unique solution $x_p^*$ in this case, and we finish the proof.
\qed
\end{lema} 

For simplicity of description, in Lemma \ref{monotonic function1}, only $x_{P_i}^0\leq x_{E_j}^0$ is considered. As for $x_{P_i}^0>x_{E_j}^0$, the similar conclusion can be obtained.

\subsection{Base Curves}
In the construction of the barrier in Section \ref{pursuitcoalition}, the following base curves are utilized, which are very essential to characterize the barrier. We emphasize that the parameters $\bm{h}_1,\bm{h}_2$ and $\bm{h}_3$ defined below represent the initial positions of three different pursuers. 

\begin{defi}[Base Curves]\label{base}\rm
For $\bm{h}_1=(x_1,y_1),\bm{h}_2=(x_2,y_2)$ and $\bm{h}_3=(x_3,y_3)$ satisfying $y_i\leq0(i=1,2,3)$ and $x_1<x_2<x_3$, define the following curves.\\
(\romannumeral1). One pursuer case:
\begin{equation}\begin{aligned}\label{base11}
&F^1_1(\bm{h}_1)=\big\{\mathbf{z}=(x,y)\in\mathbb{R}^2|\|\mathbf{z}-\bm{m}\|_2\\
&\qquad\qquad\qquad\qquad\ =\alpha\|\bm{h}_1-\bm{m}\|_2,x\leq k_1,y<0\big\} \\
&F^1_2(\bm{h}_1)=\big\{\mathbf{z}=(x,y)\in\mathbb{R}^2|(x-x_1)^2+\\
&(1-1/\alpha^2)y^2+(1-\alpha^2)y_1^2=0,x\in(k_1,k_2),y<0\big\} \\
&F^1_3(\bm{h}_1)=\big\{\mathbf{z}=(x,y)\in\mathbb{R}^2|\|\mathbf{z}-\bm{n}\|_2\\
&\qquad\qquad\qquad\qquad\ \, \ =\alpha \|\bm{h}_1-\bm{n}\|_2,x\ge k_2,y<0\big\}
\end{aligned} \end{equation}
where $k_1=\alpha^2x_1$ and  $k_2=(1-\alpha^2)l+\alpha^2 x_1$.\\
(\romannumeral2). Two pursuers case: 
\begin{equation}\begin{aligned}\label{basecurves22}
&F^2_1(\bm{h}_1,\bm{h}_2)=F_1^1(\bm{h}_1)\\
&F^2_2(\bm{h}_1,\bm{h}_2)=\big\{\mathbf{z}=(x,y)\in\mathbb{R}^2|(x-x_1)^2+\\
&\ (1-1/\alpha^2)y^2+(1-\alpha^2)y_1^2=0,x\in(k_3,k_4),y<0\big\} \\
&F^2_3(\bm{h}_1,\bm{h}_2)=\big\{\mathbf{z}=(x,y)\in\mathbb{R}^2|\|\mathbf{z}-\bm{p}_c\|_2\\
&\qquad\qquad\qquad\quad =\alpha\|\bm{h}_1-\bm{p}_c\|_2,x\in[k_4,k_5],y<0\big\}\\
&F^2_4(\bm{h}_1,\bm{h}_2)=\big\{\mathbf{z}=(x,y)\in\mathbb{R}^2|(x-x_2)^2+\\
&\ (1-1/\alpha^2)y^2+(1-\alpha^2)y_2^2=0,x\in(k_5,k_6),y<0\big\}\\
&F^2_5(\bm{h}_1,\bm{h}_2)=F_3^1(\bm{h}_2)\\
\end{aligned}\end{equation}
where $\bm{p}_c=(x_c,0)$ is the unique point on the $x$-axis such that  $\|\bm{p}_c-\bm{h}_1\|_2=\|\bm{p}_c-\bm{h}_2\|_2$. Here, $k_3=\alpha^2x_1,k_4=(1-\alpha^2)x_{c}+\alpha^2x_1,k_5=(1-\alpha^2)x_{c}+\alpha^2x_2$ and $k_6=(1-\alpha^2)l+\alpha^2x_2$.\\
(\romannumeral3). Three pursuers case:
\begin{equation}\begin{aligned}\label{basecurvethree}
&F^3(\bm{h}_1,\bm{h}_2,\bm{h}_3)=\big\{\mathbf{z}=(x,y)\in\mathbb{R}^2|(x-x_2)^2+\\
&\ (1-1/\alpha^2)y^2+(1-\alpha^2)y_2^2=0,x\in(k_7,k_8),y<0\big\}
\end{aligned}\end{equation}
where $k_7=(1-\alpha^2)x_{c1}+\alpha^2x_2$ and $k_8=(1-\alpha^2)x_{c2}+\alpha^2x_2$. Two points $\bm{p}_{c1}=(x_{c1},0)$ and $\bm{p}_{c2}=(x_{c2},0)$ are respectively given by $\|\bm{h}_1-\bm{p}_{c1}\|_2=\|\bm{h}_2-\bm{p}_{c1}\|_2$ and $\|\bm{h}_2-\bm{p}_{c2}\|_2=\|\bm{h}_3-\bm{p}_{c2}\|_2$. Although it is called three pursuers case, $F^3$ only depends on $\bm{h}_2$ while the roles of $\bm{h}_1$ and $\bm{h}_3$ are to decide two boundaries $k_7$ and $k_8$ for $x$. Thus, each point on $F^3$ only depends on at most two pursuers by cutting $F^3$ into two parts.
\end{defi}

\begin{rek}\label{rek32}\rm
It can be verified that given $\bm{h}_1,\bm{h}_2$ and $\bm{h}_3$, all base curves are explicit and smooth if they are not empty. Actually, the following discussion will show that the barrier, a focus in our paper, is composed by these curves. How to use these base curves will be stated clearly in the next section.
\end{rek}

\section {One Pursuit Coalition Versus One Evader}\label{pursuitcoalition}
\subsection{Problem Formulation}\label{noProblem Formulation}
Before investigating the number of the evaders which would be captured in $\Omega_{\rm play}$ before entering $\Omega_{\rm tar}$, a subgame between a pursuit coalition $k$ and $E_j$ is analyzed as a building block. Denote the joint game domain for the pursuit coalition $k$ by $\Omega^{n_k}=\Omega\times\cdots\times\Omega$. Similarly, denote the joint play region, target region, target line and control constraint for the pursuit coalition $k$ by $\Omega_{\rm play}^{n_k},\Omega_{\rm tar}^{n_k},\mathcal{T}^{n_k}$ and $\mathcal{U}^{n_k}$, respectively. Denote the joint control input and initial state of the pursuit coalition $k$ by $\mathcal{P}_k=\big\{\mathbf{u}_{P_{m_1}},...,\mathbf{u}_{P_{m_{n_k}}}\big\}\in\mathcal{U}^{n_k}$ and $\mathcal{X}_k^0=\big\{\mathbf{x}_{P_{m_1}}^{0},...,\mathbf{x}_{P_{m_{n_k}}}^{0}\big\}\in\mathbb{R}^{2n_k}$, respectively. 

In order to develop a basic approach, we first restrict $\mathcal{X}^0_k\in\Omega_{\rm play}^{n_k}\cup\mathcal{T}^{n_k}$ in Sections \ref{section2v1} and \ref{sectionmultiplev1}, that is, all pursuers initially lie in the play region $\Omega_{\rm play}$ or target line $\mathcal{T}$ as Assumption 2 states. Then, the case $\mathcal{X}^0_k\in\Omega^{n_k}$, that is, all pursuers can start the game from any positions in $\Omega$ as Assumption 3 states, will be discussed in Section \ref{sectionextension}. For clarity of the description, $k$ refers to the pursuit coalition $k$.

Given $\Omega$ and $\mathcal{T}$, namely, specifying $\Omega_{\rm tar}$ and $\Omega_{\rm play}$, the following problems will be addressed in this section.

\begin{pbm}\label{p31}\rm
Consider $k$ and $E_j$. Given $\mathcal{X}_k^0$, find the region $\mathcal{W}_P^{n_k}(\mathcal{X}^0_k)\subset\Omega_{\rm play}$ in which if $E_j$ initially lies, there exists a joint pursuit control input $\mathcal{P}_k\in\mathcal{U}^{n_k}$ such that $E_j$ can be captured before entering $\Omega_{\rm tar}$ regardless of its evasion control input $\mathbf{u}_{E_j}\in\mathcal{U}$.
\end{pbm}

\begin{pbm}\label{p32}\rm
Consider $k$ and $E_j$. Given $\mathcal{X}_k^0$, find the region $\mathcal{W}_E^{n_k}(\mathcal{X}^0_k)\subset\Omega_{\rm play}$ in which if $E_j$ initially lies, an evasion control input $\mathbf{u}_{E_j}\in\mathcal{U}$ exists such that $E_j$ can enter $\Omega_{\rm tar}$ without being captured regardless of the joint pursuit control $\mathcal{P}_k\in\mathcal{U}^{n_k}$.
\end{pbm}

Thus, $\mathcal{W}_P^{n_k}(\mathcal{X}^0_k)$ and $\mathcal{W}_E^{n_k}(\mathcal{X}^0_k)$ are the respective winning regions for $k$ and $E_j$, and we call them PWR and EWR respectively. According to Isaacs' book \cite{Is1967Diff}, the barrier, denoted by $\mathcal{B}^{n_k}(\mathcal{X}_k^0)$, is the surface separating $\mathcal{W}_P^{n_k}(\mathcal{X}^0_k)$ from $\mathcal{W}_E^{n_k}(\mathcal{X}^0_k)$, on which no team can guarantee its own winning. In this case, $\mathcal{B}^{n_k}(\mathcal{X}_k^0)$ is a curve, and $\mathcal{B}^{n_k}(\mathcal{X}_k^0)\cup\mathcal{W}_P^{n_k}(\mathcal{X}_k^0)\cup\mathcal{W}_E^{n_k}(\mathcal{X}_k^0)=\Omega_{\rm play}$. Unless needed for clarity, the initial condition occurring in the expressions of the barrier and winning regions will be dropped from now on.

More visually, given $\mathcal{X}_k^0$, the winning region $\mathcal{W}_P^{n_k}$ can be interpreted as the capturable region of $k$ when facing one evader, while the winning region $\mathcal{W}_E^{n_k}$ corresponds to its uncapturable region. Note that if $\mathcal{B}^{n_k}$ can be obtained, $\mathcal{W}_P^{n_k}$ and $\mathcal{W}_E^{n_k}$ split by $\mathcal{B}^{n_k}$ will come out immediately, where the region closer to $\mathcal{T}$ is $\mathcal{W}_E^{n_k}$, and the other is $\mathcal{W}_P^{n_k}$. Thus, the primary focus in this section is to construct the barrier $\mathcal{B}^{n_k}$.

For a clear symbol description, $\mathcal{B}^{n_k}(\mathcal{X}_k^0)$ is the barrier determined by $k$ with the initial position $\mathcal{X}_k^0$ and $n_k$ pursuers. More generally, $\mathcal{B}^{n_k-1}(\mathcal{X}_k^0\setminus\mathbf{x}_{P_i}^0)$ denotes the barrier determined by a pursuit coalition with the initial position $\mathcal{X}_k^0\setminus\mathbf{x}_{P_i}^0$ and $n_k-1$ pursuers, where $\mathcal{X}_k^0\setminus\mathbf{x}_{P_i}^0$ denotes the remainder in $\mathcal{X}_k^0$ when $\mathbf{x}_{P_i}^0$ is eliminated. The similar notations are also applied for $\mathcal{W}_P^{n_k}$ and $\mathcal{W}_E^{n_k}$. 

Let $\breve{\mathcal{B}}^{n_k},\breve{\mathcal{W}}_P^{n_k}$ and $\breve{\mathcal{W}}_E^{n_k}$ respectively denote the barrier, PWR, and EWR of $k$ when the boundary of $\Omega_{\rm play}$ is ignored, which will be used to construct $\mathcal{B}^{n_k},\mathcal{W}_P^{n_k}$ and $\mathcal{W}_E^{n_k}$. These three new notations have the similar meanings, citation ways and notation expressions as $\mathcal{B}^{n_k},\mathcal{W}_P^{n_k}$ and $\mathcal{W}_E^{n_k}$ respectively, and the only difference is that they are used without considering the boundary of $\Omega_{\rm play}$. They are introduced only for a clear proof.

To clarify clearly, all pursuers in $k$ are classified into two categories from the perspective of these pursuers' relative positions, which plays a crucial role in barrier construction.

\begin{defi}[Active and Inactive Pursuers]\label{active}\rm
For a pursuer $P_i(i\in\mathcal{I}_k)$ in $k$, if there exists at least one point in $\mathcal{T}$ that $P_i$ can reach before all other pursuers in $k$, we call $P_i$ is an active pursuer in $k$. Otherwise, we call $P_i$ is an inactive pursuer in $k$.
\end{defi}


\begin{rek}\label{rek31}\rm
As illustrated below, whether a pursuer is active or inactive in a pursuit coalition strictly depends on if it contributes to the barrier construction of this pursuit coalition. In brief, the contribution means being active. It can be noted that a pursuer who is inactive in one pursuit coalition may be active in another pursuit coalition, and vice versa.
\end{rek}


Next, introduce a critical payoff function which is employed to construct the barrier and provides strategies for the players. 
\begin{defi}[Payoff Function]\label{payoffADR}\rm
For $k$ and $E_j$, if $E_j$ can succeed to reach $\mathcal{T}$, take the distance of $E_j$ to the closest pursuer exactly when $E_j$ arrives at $\mathcal{T}$ as the payoff function. This payoff function $J$ and the associated value function $V$ are respectively given by
\begin{equation}\begin{aligned}\label{valueADR1v1}
&J=\min_{i\in\mathcal{I}_k}\|\mathbf{x}_{P_{i}}(t_1)-\mathbf{x}_{E_j}(t_1)\|_2,V=\min_{\mathcal{P}_k\in\mathcal{U}^{n_k}}\max_{\mathbf{u}_{E_j}\in\mathcal{U}}J
\end{aligned}\end{equation}
where $t_1$ is the first arrival time when $E_j$ reaches $\mathcal{T}$.
\end{defi}
The above payoff function, which is also called safe distance on arrival, can be interpreted as $E_j$ desires to reach $\mathcal{T}$ under the safest condition, while $k$ wants to approach $E_j$ as close as possible although the capture cannot be guaranteed. 

\subsection{Two Pursuers Versus One Evader}\label{section2v1}
We begin the discussion of the barrier construction by focusing on an important class of RA subgames with two pursuers and one evader, namely, when $k$ contains two pursuers $P_{i1}$ and $P_{i2}$. Focusing on this special case enables us to develop a scalable and analytical barrier, while also providing key insights into the barrier construction for general pursuit coalitions. 

When $\Omega_{\rm play}$ is square, the barrier $\mathcal{B}^2$ can be computed by the method proposed in \cite{Yan2017Rea}. However, the shape of $\Omega_{\rm play}$ in our current case is convex and not only restricted to square, which is beyond the scope of the previous work while quite common in general RA games. Our current work also allows the pursuers to start the game from $\Omega_{\rm tar}$ as Section \ref{sectionextension} shows, which is not involved in \cite{Yan2017Rea}. 

Next, we present a necessary condition that the optimal trajectories should satisfy if the players adopt their optimal strategies from (\ref{valueADR1v1}).
\begin{lema}[Optimal Trajectories]\label{optstraight1v1}\rm
For ${P_i}$ and $E_j$, consider the payoff function (\ref{valueADR1v1}) when only $\mathbf{x}_{P_i}^0(i\in\mathcal{I}_k)$ is considered in $\mathcal{X}_k^0$. Then, the optimal trajectories for $P_i$ and $E_j$ are both straight lines.

\emph{Proof:} The Hamiltonian function for this problem is $H=v_P\mathbf{u}_{P_{i}}\lambda_1^\mathsf{T}+v_E\mathbf{u}_{E_j}\lambda_2^\mathsf{T}$, where $\lambda_1\in\mathbb{R}^2$ and $\lambda_2\in\mathbb{R}^2$ are costate vectors. Therefore, based on the classical Issacs's method, the optimal controls $\mathbf{u}_{P_i}^*$ and $\mathbf{u}_{E_j}^*$ satisfy\begin{equation}\label{1v1Hamiltonnece}
\mathbf{u}_{P_i}^*= -\frac{\lambda_1}{\|\lambda_1\|_2},\mathbf{u}_{E_j}^* = \frac{\lambda_2}{\|\lambda_2\|_2},\dot{\lambda}_1=0,\dot{\lambda}_2=0.
\end{equation}

Thus, the optimal controls $\mathbf{u}_{P_i}^*$ and $\mathbf{u}_{E_j}^*$ are time-invariant, and their optimal trajectories are straight lines. \qed
\end{lema}

First, consider the case when $k$ contains only one pursuer, i.e., $\mathcal{X}_k^0=\mathbf{x}^0_{P_i}(1\leq i\leq N_p)$.

\begin{figure}\centering
\graphicspath{{Figures/}}
\includegraphics[width=70mm,height=48mm]{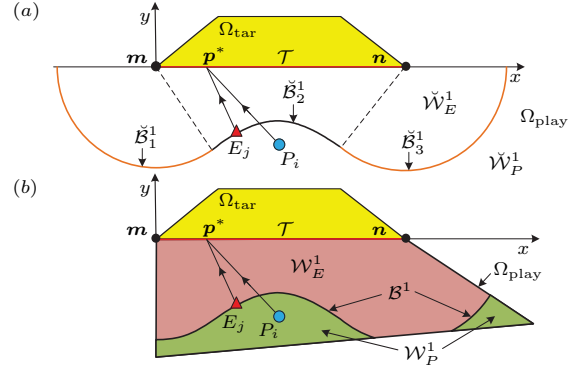}
\put(-107,9){$\scriptstyle{P_i}$}
\put(-121,12){$\scriptstyle{E_j}$}
\put(-99,73.5){$\scriptstyle{P_i}$}
\put(-119,78){$\scriptstyle{E_j}$}
\put(-52,0){$\scriptstyle{\mathcal{W}_P^1}$}
\put(-95,33){$\scriptstyle{\mathcal{W}_E^1}$}
\put(-58,23){$\scriptstyle{\mathcal{B}^1}$}
\put(-45,95){$\scriptstyle{\breve{\mathcal{W}}_E^1}$}
\put(-20,72){$\scriptstyle{\breve{\mathcal{W}}_P^1}$}
\put(-154,82){$\scriptstyle{\breve{\mathcal{B}}_1^1}$}
\put(-98,99){$\scriptstyle{\breve{\mathcal{B}}_2^1}$}
\put(-53,80){$\scriptstyle{\breve{\mathcal{B}}_3^1}$}
\put(-18,32){$\scriptstyle{\Omega_{\rm play}}$}
\put(-8,90){$\scriptstyle{\Omega_{\rm play}}$}
\put(-122,58){$\scriptstyle{\Omega_{\rm tar}}$}
\put(-122,123){$\scriptstyle{\Omega_{\rm tar}}$}
\put(-128,48){$\scriptstyle{\bm{p}^*}$}
\put(-128,113){$\scriptstyle{\bm{p}^*}$}
\put(-100,47){$\scriptstyle{\mathcal{T}}$}
\put(-100,112){$\scriptstyle{\mathcal{T}}$}
\put(-157,47){$\scriptstyle{\bm{m}}$}
\put(-157,112){$\scriptstyle{\bm{m}}$}
\put(-64,47){$\scriptstyle{\bm{n}}$}
\put(-64,112){$\scriptstyle{\bm{n}}$}
\put(-9,40){$\scriptstyle{x}$}
\put(-12,105){$\scriptstyle{x}$}
\put(-152,63){$\scriptstyle{y}$}
\put(-152,128){$\scriptstyle{y}$}
\put(-200,130){$\scriptstyle{(a)}$}
\put(-200,63){$\scriptstyle{(b)}$}
\caption{The barrier and two winning regions determined by $P_i$ when playing with $E_j$. $(a)$ Without considering the boundary of $\Omega_{\rm play}$, the barrier $\breve{\mathcal{B}}^1$ consists of three curves: $\breve{\mathcal{B}}^1_1$ (orange), $\breve{\mathcal{B}}^1_2$ (black) and $\breve{\mathcal{B}}^1_3$ (orange). The barrier $\breve{\mathcal{B}}^1$ splits $\Omega_{\rm play}$ into two winning regions $\breve{\mathcal{W}}_P^1$ and $\breve{\mathcal{W}}_E^1$. If $E_j$ lies in $\breve{\mathcal{W}}_P^1$, $P_i$ can capture $E_j$ before the latter reaches $\mathcal{T}$, while if $E_j$ lies in $\breve{\mathcal{W}}_E^1$, he can guarantee his arrival in $\mathcal{T}$. If $E_j$ lies in $\breve{\mathcal{B}}^1$, under two players' optimal strategies, $E_j$ is captured by $P_i$ exactly when reaching $\mathcal{T}$. $(b)$ By considering the boundary of $\Omega_{\rm play}$, $\mathcal{B}^1$ is given by $\breve{\mathcal{B}}^1\cap\Omega_{\rm play}$, as the optimal trajectories for $P_i$ and $E_j$ are both straight lines and $\Omega_{\rm play}$ is convex. Naturally, $\breve{\mathcal{W}}_P^1$ becomes the green region $\mathcal{W}_P^1$, and $\breve{\mathcal{W}}_E^1$ becomes the red region $\mathcal{W}_E^1$.\label{figT:3}}
\end{figure}

\begin{thom}[Barrier for One Pursuer]\label{barrierone}\rm
Consider the system (1) satisfying Assumptions \ref{isolateinitial}, \ref{constrainedinitial} and \ref{speedratio}. If a pursuit coalition $k$ only contains one pursuer $P_i$, then $\mathcal{B}^1(\mathcal{X}_k^0)=\breve{\mathcal{B}}^1\cap\Omega_{\rm play}$, where $\breve{\mathcal{B}}^1=\cup_{s=1}^3\breve{\mathcal{B}}_s^1$ and $\breve{\mathcal{B}}_s^1=F_s^1(\mathbf{x}^0_{P_i})$ for all $s=1,2,3$.
\end{thom}
\emph{Proof:} Since $k$ only contains one pursuer $P_i$, according to \autoref{binary coalitions}, we have $\mathcal{X}_k^0=\mathbf{x}^0_{P_i}$, $n_k=1$, and $\mathcal{I}_k=m_1=i$. \autoref{constrainedinitial} confines $\mathbf{x}^0_{P_i}$ in $\Omega_{\rm play}\cup\mathcal{T}$. First, we do not consider the boundary of $\Omega_{\rm play}$, namely, focus on the computation of $\breve{\mathcal{B}}^1$ defined in Section \ref{noProblem Formulation} which is proved to consist of three parts $\breve{\mathcal{B}}^1_1,\breve{\mathcal{B}}^1_2$ and $\breve{\mathcal{B}}^1_3$ in the following discussion, as depicted in Fig. \ref{figT:3}(a). 

Assume that $E_j$ can succeed to reach $\mathcal{T}$, and denote $E_j$'s optimal target point (OTP) in $\mathcal{T}$ by $\bm{p}^*=(x_p^*,0)$ such that $J$ in (\ref{valueADR1v1}) is maximized. It can be seen that in this case, (\ref{valueADR1v1}) is simplified to
\begin{equation}\begin{aligned}\label{onepursuerpayoffAWR}
&J=\|\mathbf{x}_{P_i}(t_1)-\mathbf{x}_{E_j}(t_1)\|_2,V=\min_{\mathbf{u}_{P_i}\in\mathcal{U}}\max_{\mathbf{u}_{E_j}\in\mathcal{U}}J
\end{aligned}\end{equation}
where $t_1$ is the first arrival time when $E_j$ reaches $\mathcal{T}$.

It follows from \autoref{optstraight1v1} that for the payoff function (\ref{onepursuerpayoffAWR}), the optimal trajectories for $P_i$ and $E_j$ are both straight lines. Also note that the non-anticipative information structure stated in Section \ref{assumtionsection} implies that $P_i$ knows $E_j$'s current input (i.e., direction of motion). Since $P_i$ aims to minimize $J$ in (\ref{onepursuerpayoffAWR}), the optimal strategy for $P_i$ is to move towards the same target point in $\mathcal{T}$ as $E_j$ does, as Fig. \ref{figT:3}(a) illustrates. Hence, if $E_j$ selects $\bm{p}$ in $\mathcal{T}$ to go, the payoff function $J$ in (\ref{onepursuerpayoffAWR}) is equivalent to the function\begin{equation}\label{function1v1}
G_1(x_p)=\|\bm{p}-\mathbf{x}_{P_i}^0\|_2-\frac{\|\bm{p}-\mathbf{x}_{E_j}^0\|_2}{\alpha},\bm{p}=(x_p,0)\in\mathcal{T}
\end{equation}
which is $P_i$'s distance to $E_j$ exactly when $E_j$ arrives at $\bm{p}$, if $P_i$ and $E_j$ both move directly towards $\bm{p}$. Notice that $0\leq x_p\leq l$, where $l$ is the length of $\mathcal{T}$.

Note that $E_j$ strives to maximize $J$, i.e., $G_1(x_p)$, as (\ref{onepursuerpayoffAWR}) shows. Thus, if $P_i$ and $E_j$ both adopt their optimal strategies, $x_p^*$ must be a maximum point of $G_1(x_p)$ for $x_p\in[0,l]$. Note that the monotony of $G_1(x_p)$ has been shown in Lemma \ref{monotonic function1}. If the maximum point $x_p^*$ of $G_1(x_p)$ for $x_p\in[0,l]$ occurs in the interior of $[0,l]$, i.e., $(0,l)$, then $x_p^*$ is an extreme point of $G_1(x_p)$, that is, 
\begin{equation}\label{FNC}
\frac{{\rm d}G_1(x_p^*)}{{\rm d}x_p}=0\Rightarrow\frac{x_p^*-x_{P_i}^0}{\|\bm{p}^*-\mathbf{x}_{P_i}^0\|_2}=\frac{x_p^*-x_{E_j}^0}{\alpha\|\bm{p}^*-\mathbf{x}_{E_j}^0\|_2}.
\end{equation}
\hspace{0.2 in}Furthermore, assume $\mathbf{x}_{E_j}^0\in\breve{\mathcal{B}}^1$, as Fig. \ref{figT:3}(a) shows. Then, $E_j$ will be captured by $P_i$ exactly when reaching $\mathcal{T}$ under two players' optimal strategies from (\ref{onepursuerpayoffAWR}), namely, $V=0$. Thus, $P_i$ and $E_j$ will reach $\bm{p}^*$ simultaneously as follows:
\begin{equation}\label{distanceunequ}
G_1(x_p^*)=0\Rightarrow\|\bm{p}^*-\mathbf{x}_{E_j}^0\|_2=\alpha\|\bm{p}^*-\mathbf{x}_{P_i}^0\|_2.
\end{equation}
This feature (\ref{distanceunequ}) reflects that when $\mathbf{x}_{E_j}^0\in\breve{\mathcal{B}}^1$, under two players' optimal strategies from (\ref{onepursuerpayoffAWR}), no player can guarantee its winning, namely, the capture and arrival happen at the same time. 

The following analysis will separately consider three cases: $x_p^*\in(0,l),x_p^*=0$ and $x_p^*=l$, due to their different features.

First, consider $x_p^*\in(0,l)$. It has been shown above that in this case, $x_p^*$ is an extreme point of $G_1(x_p)$. Denote this part of $\breve{\mathcal{B}}^1$ by $\breve{\mathcal{B}}^1_2$ which is in black shown in Fig. \ref{figT:3}(a). Next, we compute the expression of $\breve{\mathcal{B}}^1_2$.

Substituting (\ref{distanceunequ}) into (\ref{FNC}) yields\begin{equation}\label{OTPunequ11}
\alpha^2(x_p^*-x_{P_i}^0)=x_p^*-x_{E_j}^0\Rightarrow x_p^*=\frac{x_{E_j}^0-\alpha^2x_{P_i}^0}{1-\alpha^2}.
\end{equation}
Then, substituting $x_p^*$ given by (\ref{OTPunequ11}) into (\ref{distanceunequ}) leads to
\begin{equation}\begin{aligned}\label{barrier1v1explicit}
&\Big\|\big(\frac{x_{E_j}^0-\alpha^2x_{P_i}^0}{1-\alpha^2},0\big)-\mathbf{x}_{E_j}^0\Big\|_2^2\\
&\qquad\qquad\qquad\qquad=\alpha^2\Big\|\big(\frac{x_{E_j}^0-\alpha^2x_{P_i}^0}{1-\alpha^2},0\big)-\mathbf{x}_{P_i}^0\Big\|_2^2\Rightarrow\\
&(x_{E_j}^0-x_{P_i}^0)^2+(1-1/\alpha^2)(y_{E_j}^0)^2+(1-\alpha^2)(y_{P_i}^0)^2 =0
\end{aligned}\end{equation}
which characterizes the relationship between the initial positions of $P_i$ and $E_j$ when $\mathbf{x}_{E_j}^0\in\breve{\mathcal{B}}^1_2$. Since $x_p^*\in(0,l)$, (\ref{OTPunequ11}) implies 
\begin{equation}\begin{aligned}\label{barrier1v1explicitinterral}
x_{E_j}^0\in\big(\alpha^2x_{P_i}^0,(1-\alpha^2)l+\alpha^2x_{P_i}^0\big).
\end{aligned}\end{equation}

Also note that \autoref{constrainedinitial} confines $\mathbf{x}^0_{E_j}$ in $\Omega_{\rm play}$, that is, $y_{E_j}^0<0$. Hence, given $\mathbf{x}_{P_i}^0$, by taking all positions $\mathbf{x}_{E_j}^0$ satisfying (\ref{barrier1v1explicit}) and (\ref{barrier1v1explicitinterral}) with $y_{E_j}^0<0$, all initial positions of $E_j$ lying in $\breve{\mathcal{B}}^1_2$ are found. Equivalently, these initial positions of $E_j$ form the $\breve{\mathcal{B}}^1_2$. Thus, by replacing $\mathbf{x}_{E_j}^0$ with a general variable $\mathbf{z}\in\mathbb{R}^2$, it follows from (\ref{barrier1v1explicit}), (\ref{barrier1v1explicitinterral}) and $y_{E_j}^0<0$ that $\breve{\mathcal{B}}_2^1$ is as follows:
\begin{equation}\begin{aligned}\label{barrier1v1expression11}
\breve{\mathcal{B}}_2^1=\big\{\mathbf{z}=(x,y)\in\mathbb{R}^2|(x-x_{P_i}^0)^2+(1-1/\alpha^2)y^2\\
+(1-\alpha^2)(y_{P_i}^0)^2 =0,x\in(k_1,k_2),y<0\big\}
\end{aligned}\end{equation}
where $k_1=\alpha^2x_{P_i}^0$ and $k_2=(1-\alpha^2)l+\alpha^2x_{P_i}^0$. For clarity, it can be seen that $\breve{\mathcal{B}}_2^1$ can be expressed by a base curve in (\ref{base11}), that is, $\breve{\mathcal{B}}_2^1=F_2^1(\mathbf{x}^0_{P_i})$.

Then, we focus on the case $x_p^*=0$, namely, $\bm{p}^*=\bm{m}$. Denote this part of $\breve{\mathcal{B}}^1$ by $\breve{\mathcal{B}}^1_1$ which is the left orange curve in Fig. \ref{figT:3}(a). Next, we compute the expression of $\breve{\mathcal{B}}^1_1$.

In this scenario, (\ref{distanceunequ}) still holds and becomes
 \begin{equation}\begin{aligned}\label{barrier1v1expression12233}
\|\mathbf{x}_{E_j}^0-\bm{m}\|_2=\alpha\|\mathbf{x}_{P_i}^0-\bm{m}\|_2.
 \end{aligned}\end{equation}
Note that $x_{E_j}^0\leq\alpha^2x_{P_{i}}^0$ is straightforward based on the interval (\ref{barrier1v1explicitinterral}) of $x_{E_j}^0$ for $\mathbf{x}_{E_j}^0\in\breve{\mathcal{B}}^1_2$ and the monotony of $G_1(x_p)$ given by \autoref{monotonic function1}, as Fig. \ref{figT:3}(a) shows. Thus, by replacing $\mathbf{x}_{E_j}^0$ with a general variable $\mathbf{z}\in\mathbb{R}^2$, it follows from (\ref{barrier1v1expression12233}) and $x_{E_j}^0\leq\alpha^2x_{P_{i}}^0$ that $\breve{\mathcal{B}}_1^1$ is as follows:\begin{equation}\begin{aligned}\label{barrier1v1b123}
 \breve{\mathcal{B}}_1^1=\big\{&\mathbf{z}=(x,y)\in\mathbb{R}^2|\|\mathbf{z}-\bm{m}\|_2\\
 &=\alpha\|\mathbf{x}_{P_i}^0-\bm{m}\|_2,x\leq k_1,y<0\big\}
\end{aligned}\end{equation} 
which can also be represented by a base curve in (\ref{base11}), that is, $\breve{\mathcal{B}}_1^1=F_1^1(\mathbf{x}^0_{P_i})$.

If $x_p^*=l$, namely, $\bm{p}^*=\bm{n}$, denote this part of $\breve{\mathcal{B}}^1$ by $\breve{\mathcal{B}}^1_3$ which is the right orange curve in Fig. \ref{figT:3}(a). Similar to the case $x_p^*=0$, $\breve{\mathcal{B}}_3^1$ is as follows:\begin{equation}\begin{aligned}\label{barrier1v1b12333}
 \breve{\mathcal{B}}_3^1=\big\{&\mathbf{z}=(x,y)\in\mathbb{R}^2|\|\mathbf{z}-\bm{n}\|_2\\
 &=\alpha\|\mathbf{x}_{P_i}^0-\bm{n}\|_2,x\ge k_2,y<0\big\}.
\end{aligned}\end{equation} 
Thus, $\breve{\mathcal{B}}_3^1=F_3^1(\mathbf{x}^0_{P_i})$ can be obtained.

Therefore, we have constructed the barrier $\breve{\mathcal{B}}^1=\cup_{s=1}^3\breve{\mathcal{B}}_s^1$, without considering the boundary of $\Omega_{\rm play}$, shown in Fig. \ref{figT:3}(a). Then, consider the effect of the boundary of $\Omega_{\rm play}$. Since  the optimal trajectories for two players are both straight lines and $\Omega_{\rm play}$ is convex, it can be concluded that $\mathcal{B}^1=\breve{\mathcal{B}}^1\cap\Omega_{\rm play}$, as Fig. \ref{figT:3}(b) shows. The shape of $\Omega_{\rm play}$ in Fig. \ref{figT:3}(b), which is general, is introduced just for showing that this theorem is applied for any convex $\Omega_{\rm play}$. 

As depicted in Fig. \ref{figT:3}(b), since $\mathcal{W}_P^1$ and $\mathcal{W}_E^1$ are two subregions of $\Omega_{\rm play}$ and separated by $\mathcal{B}^1$, and $\mathcal{W}_E^1$ is closer to $\mathcal{T}$ than $\mathcal{W}_P^1$, then the green region is the PWR $\mathcal{W}_P^1$ and the red region is the EWR $\mathcal{W}_E^1$. Thus, we finish the proof. \qed



Next, we consider the case when $k$ contains two pursuers, i.e., $\mathcal{X}_k^0=\big\{\mathbf{x}^0_{P_{i1}},\mathbf{x}^0_{P_{i2}}\big\}(1\leq i1,i2\leq N_p,i1\neq i2)$. The main result of this section is presented below, which gives the analytical barrier for two pursuers case.

\begin{thom}[Barrier for Two Pursuers]\label{barrier two}\rm
Consider the system (1) and suppose that Assumptions \ref{isolateinitial}, \ref{constrainedinitial} and \ref{speedratio} hold. If a pursuit coalition $k$ contains two pursuers $P_{i1}$ and $P_{i2}$, its barrier can be analytically calculated as follows:
\begin{itemize}
\item[a] (Only one active pursuer). If only ${P_i}$ is an active pursuer, $\mathcal{B}^2(\mathcal{X}_k^0)=\mathcal{B}^1(\mathbf{x}^0_{P_i})$ where $i=i1$ or $i2$. \\
\item[b] (Two active pursuers). If $P_{i1}$ and ${P_{i2}}$ are both active pursuers, assume $x^0_{P_{i1}}<x^0_{P_{i2}}$. Then, $\mathcal{B}^2(\mathcal{X}_k^0)=\breve{\mathcal{B}}^2\cap\Omega_{\rm play}$, where $\breve{\mathcal{B}}^2=\cup_{s=1}^5\breve{\mathcal{B}}_s^2$ and $\breve{\mathcal{B}}_s^2=F_s^2(\mathbf{x}^0_{P_{i1}},\mathbf{x}^0_{P_{i2}})$ for all $s=1,...,5$.
\end{itemize}
\end{thom}

\begin{figure}\centering
\graphicspath{{Figures/}}
\includegraphics[width=70mm,height=52mm]{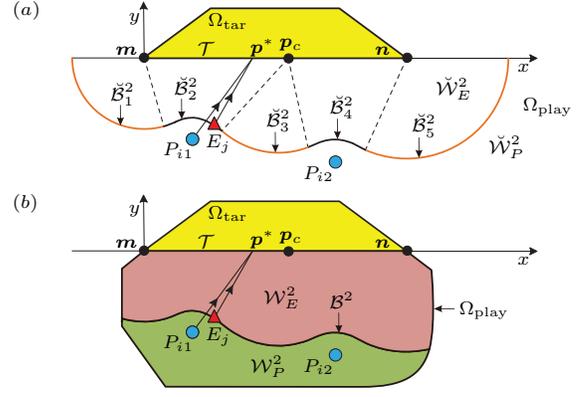}
\put(-143,16){$\scriptstyle{P_{i1}}$}
\put(-90,7){$\scriptstyle{{P_{i2}}}$}
\put(-127,19){$\scriptstyle{E_j}$}
\put(-143,89){$\scriptstyle{P_{i1}}$}
\put(-90,80){$\scriptstyle{{P_{i2}}}$}
\put(-127,92){$\scriptstyle{E_j}$}
\put(-110,7){$\scriptstyle{\mathcal{W}_P^2}$}
\put(-105,33){$\scriptstyle{\mathcal{W}_E^2}$}
\put(-80,31){$\scriptstyle{\mathcal{B}^2}$}
\put(-40,112){$\scriptstyle{\breve{\mathcal{W}}_E^2}$}
\put(-20,90){$\scriptstyle{\breve{\mathcal{W}}_P^2}$}
\put(-163,110){$\scriptstyle{\breve{\mathcal{B}}_1^2}$}
\put(-139,112){$\scriptstyle{\breve{\mathcal{B}}_2^2}$}
\put(-104,100){$\scriptstyle{\breve{\mathcal{B}}_3^2}$}
\put(-80,105){$\scriptstyle{\breve{\mathcal{B}}_4^2}$}
\put(-49,98){$\scriptstyle{\breve{\mathcal{B}}_5^2}$}
\put(-31,29){$\scriptstyle{\Omega_{\rm play}}$}
\put(-7,106){$\scriptstyle{\Omega_{\rm play}}$}
\put(-126,65){$\scriptstyle{\Omega_{\rm tar}}$}
\put(-126,138){$\scriptstyle{\Omega_{\rm tar}}$}
\put(-110,55){$\scriptstyle{\bm{p}^*}$}
\put(-110,128){$\scriptstyle{\bm{p}^*}$}
\put(-130,54){$\scriptstyle{\mathcal{T}}$}
\put(-130,127){$\scriptstyle{\mathcal{T}}$}
\put(-99,57){$\scriptstyle{\bm{p}_c}$}
\put(-99,130){$\scriptstyle{\bm{p}_c}$}
\put(-161,54){$\scriptstyle{\bm{m}}$}
\put(-161,127){$\scriptstyle{\bm{m}}$}
\put(-63,54){$\scriptstyle{\bm{n}}$}
\put(-63,127){$\scriptstyle{\bm{n}}$}
\put(-9,47.5){$\scriptstyle{x}$}
\put(-9,120.5){$\scriptstyle{x}$}
\put(-156,67){$\scriptstyle{y}$}
\put(-156,140){$\scriptstyle{y}$}
\put(-200,142){$\scriptstyle{(a)}$}
\put(-200,69){$\scriptstyle{(b)}$}
\caption{The barrier and two winning regions determined by two active pursuers $P_{i1}$ and $P_{i2}$ when playing with $E_j$, in which two pursuers both contribute to the construction of $\breve{\mathcal{B}}^2$. $(a)$ Without considering the boundary of $\Omega_{\rm play}$, the barrier $\breve{\mathcal{B}}^2$ consists of five parts and is given by $\cup_{s=1}^5\breve{\mathcal{B}}^2_s$. Two winning regions $\breve{\mathcal{W}}_P^2$ and $\breve{\mathcal{W}}_E^2$ for each team's guaranteed winning are split by $\breve{\mathcal{B}}^2$. $(b)$ By considering the boundary of $\Omega_{\rm play}$, $\mathcal{B}^2$ is given by $\breve{\mathcal{B}}^2\cap\Omega_{\rm play}$, as the optimal trajectories for the players who work in the payoff function, are straight lines and $\Omega_{\rm play}$ is convex. Naturally, $\breve{\mathcal{W}}_P^2$ becomes the green region $\mathcal{W}_P^2$, and $\breve{\mathcal{W}}_E^2$ becomes the red region $\mathcal{W}_E^2$.\label{figT:4}}
\end{figure}

\emph{Proof:} Since $k$ contains two pursuers $P_{i1}$ and $P_{i2}$, according to \autoref{binary coalitions}, we have $\mathcal{X}_k^0=\{\mathbf{x}^0_{P_{i1}},\mathbf{x}^0_{P_{i2}}\}$, $n_k=2$, and $\mathcal{I}_k=\{m_1,m_2\}=\{i1,i2\}$. Similar to the proof of \autoref{barrierone}, ignore the boundary of $\Omega_{\rm play}$ and consider $\breve{\mathcal{B}}^2$ first, as shown in Fig. \ref{figT:4}(a).

Part (a): It suffices to consider the case where $P_{i1}$ is an active pursuer and $P_{i2}$ is an inactive pursuer. 

Since $P_{i2}$ is an inactive pursuer, it follows from \autoref{active} that\begin{equation}\label{theromequ12}
\|\mathbf{x}^0_{P_{i1}}-\bm{p}\|_2\leq\|\mathbf{x}^0_{P_{i2}}-\bm{p}\|_2
\end{equation}
holds for all $\bm{p}\in\mathcal{T}$. Hence, $P_{i1}$ can reach any point in $\mathcal{T}$ no later than $P_{i2}$, and this feature guarantees that $P_{i1}$ can determine the barrier $\breve{\mathcal{B}}^2$ alone. Thus, $\breve{\mathcal{B}}^2=\breve{\mathcal{B}}^1(\mathbf{x}_{P_{i1}}^0)$.

Part (b): For two active pursuers case, we will show that $\breve{\mathcal{B}}^2$ includes five parts $\breve{\mathcal{B}}^2_s(s=1,...,5)$, as depicted in Fig. \ref{figT:4}(a). 

Since two pursuers are both active, it follows from \autoref{active} that there must exist two points $\bm{p}_1$ and $\bm{p}_2$ in $\mathcal{T}$ such that\begin{equation}\label{theromequ1233}
\|\mathbf{x}^0_{P_{i1}}-\bm{p}_1\|_2<\|\mathbf{x}^0_{P_{i2}}-\bm{p}_1\|_2,\|\mathbf{x}^0_{P_{i2}}-\bm{p}_2\|_2<\|\mathbf{x}^0_{P_{i1}}-\bm{p}_2\|_2.
\end{equation}
The former in (\ref{theromequ1233}) implies that $\breve{\mathcal{B}}^2$ depends on $P_{i1}$, and the latter in (\ref{theromequ1233}) guarantees that $\breve{\mathcal{B}}^2$ depends on $P_{i2}$. Thus, in this case, $\breve{\mathcal{B}}^2$ depends on both two pursuers. Naturally, there exists a unique point $\bm{p}_{c}=(x_{c},0)$ in $\mathcal{T}$ that two pursuers can reach at the same time, as Fig. \ref{figT:4}(a) shows.

Note that in this part, $x^0_{P_{i1}}\neq x^0_{P_{i2}}$; otherwise, it can be easily verified that one pursuer can reach any point in $\mathcal{T}$ no later than the other, which corresponds to the Part (a). Without loss of generality, assume $x^0_{P_{i1}}<x^0_{P_{i2}}$ as this theorem states. Thus, $\bm{p}_{c}$ is given by\begin{equation}\begin{aligned}\label{pcv2v1}
\|\mathbf{x}_{P_{i1}}^0-\bm{p}_{c}\|_2=\|\mathbf{x}_{P_{i2}}^0-\bm{p}_{c}\|_2\Rightarrow x_{c}=\frac{\|\mathbf{x}_{P_{i2}}^0\|_2^2-\|\mathbf{x}_{P_{i1}}^0\|_2^2}{2(x^0_{P_{i2}}-x^0_{P_{i1}})}.
\end{aligned}\end{equation}

If $x_c=0$, i.e., $\bm{p}_c=\bm{m}$, it can be seen from Fig. \ref{figT:4}(a) that $P_{i2}$ can reach any point in $\mathcal{T}$ no later than $P_{i1}$. However, $P_{i1}$ is an active pursuer, so $x_c\neq0$. Similarly, $x_c\neq l$ can be obtained. Thus, we have $0<x_c<l$. 

Take a point $\bm{p}=(x_p,0)$ in $\mathcal{T}$. It can be observed in Fig. \ref{figT:4}(a) that if $x_p\in[0,x_c)$, $\|\mathbf{x}^0_{P_{i1}}-\bm{p}\|_2<\|\mathbf{x}^0_{P_{i2}}-\bm{p}\|_2$, and if $x_p\in(x_c,l]$, $\|\mathbf{x}^0_{P_{i2}}-\bm{p}\|_2<\|\mathbf{x}^0_{P_{i1}}-\bm{p}\|_2$.

Assume $E_j$ can succeed to reach $\mathcal{T}$, and denote $E_j$'s OTP in $\mathcal{T}$ by $\bm{p}^*=(x_p^*,0)$ such that $J$ in (\ref{valueADR1v1}) is maximized. It can be noted that in this case,  (\ref{valueADR1v1}) becomes
\begin{equation}\begin{aligned}\label{valueADR2v1barrier}
J&=\min_{i=i1,i2}\|\mathbf{x}_{P_i}(t_1)-\mathbf{x}_{E_j}(t_1)\|_2\\
V&=\min_{\mathbf{u}_{P_{i1}},\mathbf{u}_{P_{i2}}\in\mathcal{U}}\max_{\mathbf{u}_{E_j}\in\mathcal{U}}J
\end{aligned}\end{equation}
where $t_1$ is the first arrival time when $E_j$ reaches $\mathcal{T}$. Furthermore, assume $\mathbf{x}_{E_j}^0\in\breve{\mathcal{B}}^2$. Therefore, under three players' optimal strategies from (\ref{valueADR2v1barrier}), $E_j$ will be captured exactly when reaching $\mathcal{T}$, namely, $V=0$.

If $x_p^*\in[0,x_c)$, as Fig. \ref{figT:4}(a) shows, this part of $\breve{\mathcal{B}}^2$ is determined by $P_{i1}$ alone, as $P_{i1}$ can reach $\bm{p}^*$ before $P_{i2}$. Thus, by following the proof of \autoref{barrierone}, it can be obtained that if $x_p^*=0$, i.e., $\bm{p}^*=\bm{m}$, and denote the related barrier by $\breve{\mathcal{B}}^2_1$, then $\breve{\mathcal{B}}^2_1=\breve{\mathcal{B}}^1_1(\mathbf{x}_{P_{i1}}^0)$, which is the leftmost orange curve in Fig. \ref{figT:4}(a) and given by (\ref{barrier1v1b123}) with $\mathbf{x}_{P_i}^0$ replaced by $\mathbf{x}_{P_{i1}}^0$. Thus, it follows from (\ref{basecurves22}) and \autoref{barrierone} that $\breve{\mathcal{B}}^2_1$ can be expressed by a base curve, i.e., $\breve{\mathcal{B}}^2_1=F_1^2(\mathbf{x}_{P_{i1}}^0,\mathbf{x}_{P_{i2}}^0)$.

Consider the remainder $x_p^*\in(0,x_c)$ and denote the related barrier by $\breve{\mathcal{B}}^2_2$, which is the leftmost black curve in Fig. \ref{figT:4}(a). Next, we compute the expression of $\breve{\mathcal{B}}^2_2$.

Since $x_p^*\in(0,x_c)$ and $J$ in (\ref{valueADR2v1barrier}) only depends on $P_{i1}$, it follows from the proof of \autoref{barrierone} that (\ref{OTPunequ11}) and (\ref{barrier1v1explicit}) still hold for $P_{i1}$ and $E_j$, that is, $x_p^*=\frac{x_{E_j}^0-\alpha^2x_{P_{i1}}^0}{1-\alpha^2}$ and
\begin{equation}\begin{aligned}\label{barrier2v1explicitDi1}
(x_{E_j}^0-x_{P_{i1}}^0)^2+(1-1/\alpha^2)(y_{E_j}^0)^2+(1-\alpha^2)(y_{P_{i1}}^0)^2 =0
\end{aligned}\end{equation}
representing the relationship between the initial positions of $P_{i1}$ and $E_j$ when $\mathbf{x}_{E_j}^0\in\breve{\mathcal{B}}^2_2$. From $x_p^*\in(0,x_c)$, we have
\begin{equation}\begin{aligned}\label{barrier2v1intervalb22}
x_{E_j}^0\in\big(\alpha^2x_{P_{i1}}^0,(1-\alpha^2)x_c+\alpha^2x_{P_{i1}}^0\big).
\end{aligned}\end{equation}
Hence, by replacing $\mathbf{x}_{E_j}^0$ with a general variable $\mathbf{z}\in\mathbb{R}^2$, it follows from (\ref{barrier2v1explicitDi1}), (\ref{barrier2v1intervalb22}) and $y_{E_j}^0<0$ that $\breve{\mathcal{B}}_2^2$ is as follows:
\begin{equation}\begin{aligned}\label{barrier2v1B22}
\breve{\mathcal{B}}_2^2=\big\{\mathbf{z}=(x,y)\in\mathbb{R}^2|(x-x_{P_{i1}}^0)^2+(1-1/\alpha^2)y^2\\
+(1-\alpha^2)(y_{P_{i1}}^0)^2 =0,x\in(k_3,k_4),y<0\big\}
\end{aligned}\end{equation}
where $k_3=\alpha^2x_{P_{i1}}^0$ and $k_4=(1-\alpha^2)x_c+\alpha^2x_{P_{i1}}^0$. For clarity, it can be seen that $\breve{\mathcal{B}}_2^2$ can be expressed by a base curve in (\ref{basecurves22}), that is, $\breve{\mathcal{B}}^2_2=F_2^2(\mathbf{x}_{P_{i1}}^0,\mathbf{x}_{P_{i2}}^0)$.

Analogously, when $P_{i2}$ can reach $\bm{p}^*$ before $P_{i1}$, i.e., $x_p^*\in(x_c,l]$, as Fig. \ref{figT:4}(a) shows, $\breve{\mathcal{B}}^2_4=F_4^2(\mathbf{x}_{P_{i1}}^0,\mathbf{x}_{P_{i2}}^0)$ and $\breve{\mathcal{B}}^2_5=F_5^2(\mathbf{x}_{P_{i1}}^0,\mathbf{x}_{P_{i2}}^0)$ can be derived respectively for $x_p^*\in(x_c,l)$ and $x_p^*=l$, where $\breve{\mathcal{B}}^2_4$ is the rightmost black curve and $\breve{\mathcal{B}}^2_5$ is the rightmost orange curve.

Finally, consider the case $x_p^*=x_c$, i.e., $\bm{p}^*=\bm{p}_c$, which corresponds to the middle orange barrier $\breve{\mathcal{B}}^2_3$ in Fig. \ref{figT:4}(a). Note that (\ref{distanceunequ}) still holds for $P_{i1}$ and $E_j$ and can be rewritten as 
\begin{equation}\label{distanceunequ22}
\|\bm{p}_c-\mathbf{x}_{E_j}^0\|_2=\alpha\|\bm{p}_c-\mathbf{x}_{P_{i1}}^0\|_2.
\end{equation}
Naturally, $x_{E_j}^0$ should lie in $[k_4,k_5]$, where $k_5=(1-\alpha^2)x_{c}+\alpha^2x_{P_{i2}}^0$ which actually is the boundary of the $x$ value of $\breve{\mathcal{B}}^2_4$. Thus, by replacing $\mathbf{x}_{E_j}^0$ with a general variable $\mathbf{z}\in\mathbb{R}^2$, it follows from (\ref{distanceunequ22}) and $y_{E_j}^0<0$ that $\breve{\mathcal{B}}_3^2$ is as follows:
\begin{equation}\begin{aligned}\label{barrier2v1B23}
 \breve{\mathcal{B}}_3^2=\big\{&\mathbf{z}=(x,y)\in\mathbb{R}^2|\|\mathbf{z}-\bm{p}_c\|_2\\
 &=\alpha\|\mathbf{x}_{P_{i1}}^0-\bm{p}_c\|_2,x\in[k_4,k_5],y<0\big\}.
\end{aligned}\end{equation}
It can also be noted that $\breve{\mathcal{B}}_3^2$ can be expressed by a base curve in (\ref{basecurves22}), that is, $\breve{\mathcal{B}}^2_3=F_3^2(\mathbf{x}_{P_{i1}}^0,\mathbf{x}_{P_{i2}}^0)$.

Therefore, we have constructed the barrier $\breve{\mathcal{B}}^2=\cup_{s=1}^5\breve{\mathcal{B}}_s^2$, without considering the boundary of $\Omega_{\rm play}$, shown in Fig. \ref{figT:4}(a). Now, consider the effect of the boundary of $\Omega_{\rm play}$. Since  the optimal trajectories for three players (if it works in the payoff function $J$) are straight lines and $\Omega_{\rm play}$ is convex, $\mathcal{B}^2=\breve{\mathcal{B}}^2\cap\Omega_{\rm play}$ holds, as Fig. \ref{figT:4}(b) shows. Intuitively, the PWR $\mathcal{W}_P^2$ is the green region, and the EWR $\mathcal{W}_E^2$ is the red region.
\qed


\subsection{General Pursuit Coalitions Versus One Evader}\label{sectionmultiplev1}
The objective in this section is to construct the barrier for general pursuit coalitions with the inspiration from the results derived in two pursuers case. First, the barrier for a special class of pursuit coalitions whose pursuers are all active, is constructed. Then theoretically, it is demonstrated that every pursuit coalition can be degenerated into a unique pursuit subcoalition of this class. Finally, this pursuit subcoalition and the original pursuit coalition are proved to possess the same barrier. To this end, the definition of this special class pursuit coalition is formally introduced as follows.

\begin{defi}[Full-Active Pursuit Coalition]\label{full active}\rm
Given a pursuit coalition $k$, if for any pursuer ${P_i}(i\in\mathcal{I}_k)$ in $k$, there always exists a point in $\mathcal{T}$ that $P_i$ can reach before all other pursuers in $k$, then $k$ is called a full-active pursuit coalition.
\end{defi}

It can be verified that if a pursuit coalition $k$ is full-active, its members must have different $x$-coordinates. Otherwise, assume that in $k$, $P_{i1}$ and $P_{i2}$ have the same $x$-coordinate. Thus, it can be noted that one pursuer can reach any point in $\mathcal{T}$ no later than the other one, which contradicts with the fact that two pursuers are both active. Hence, for every full-active pursuit coalition $k$ with $n_k\ge2$, introduce an auxiliary index set $\mathcal{\hat{I}}_k=\{\hat{m}_j|1\leq \hat{m}_j\leq N_p,j=1,2,...,n_k\}=\mathcal{I}_k$ such that $x_{P_{\hat{m}_i}}^0<x_{P_{\hat{m}_{i+1}}}^0$ holds for all $i=1,2,...,n_k-1$. More intuitively, we mean that $P_{\hat{m}_i}$ lies at the left side of $P_{\hat{m}_{i+1}}$ along the $x$-axis. We call $\mathcal{\hat{I}}_k$ as $x$-rank index set, and obviously it is unique. 

With the above $x$-rank index set $\mathcal{\hat{I}}_k$, the theorem presented below provides a scheme to construct the barrier for full-active pursuit coalitions. 

\begin{thom}[Barrier for Full-Active Pursuit Coalition]\label{barrierfull}\rm
Consider the system (1) and suppose that Assumptions \ref{isolateinitial}, \ref{constrainedinitial} and \ref{speedratio} hold. If a pursuit coalition $k$ with the initial condition $\mathcal{X}_k^0$ and $n_k$ pursuers, is full-active, its corresponding barrier can be constructed analytically as follows:\\
(\romannumeral1). If $n_k=1$, $\mathcal{B}^{n_k}(\mathcal{X}_k^0)$ is given by \autoref{barrierone}. If $n_k=2$, $\mathcal{B}^{n_k}(\mathcal{X}_k^0)$ is given by \autoref{barrier two}(b).\\
(\romannumeral2). If $3\leq n_k\leq N_p$, let $\mathcal{\hat{I}}_k$ denote its $x$-rank index set. Then, $\mathcal{B}^{n_k}(\mathcal{X}_k^0)=\breve{\mathcal{B}}^{n_k}\cap\Omega_{\rm play}$, where
\begin{equation}\begin{aligned}\label{barriergeneralnv1}
\breve{\mathcal{B}}^{n_k}=&\cup_{i=1}^2F^2_i(\mathbf{x}^0_{P_{\hat{m}_1}},\mathbf{x}^0_{P_{\hat{m}_2}})\cup_{i=1}^{n_k-1}F^2_3(\mathbf{x}^0_{P_{\hat{m}_i}},\mathbf{x}^0_{P_{\hat{m}_{i+1}}})\\
&\cup_{i=2}^{n_k-1}F^3(\mathbf{x}^0_{P_{\hat{m}_{i-1}}},\mathbf{x}^0_{P_{\hat{m}_i}},\mathbf{x}^0_{P_{\hat{m}_{i+1}}})\\
&\cup_{i=4}^5F_i^2(\mathbf{x}^0_{P_{\hat{m}_{n_k-1}}},\mathbf{x}^0_{P_{\hat{m}_{n_k}}}).
\end{aligned}\end{equation}
\end{thom}
\emph{Proof:} Note that the conclusion (\romannumeral1) is straightforward. Thus, we focus on the conclusion (\romannumeral2). Construct $\breve{\mathcal{B}}^{n_k}$ first, namely, ignore the boundary of $\Omega_{\rm play}$. 

Since $k$ is full-active, according to \autoref{full active}, for every pursuer in $k$, there always exists a point in $\mathcal{T}$ that it can reach before all other pursuers in $k$. Thus,  every pursuer in $k$ contributes to $\breve{\mathcal{B}}^{n_k}$.

\begin{figure}\centering
\graphicspath{{Figures/}}
\includegraphics[width=70mm,height=52mm]{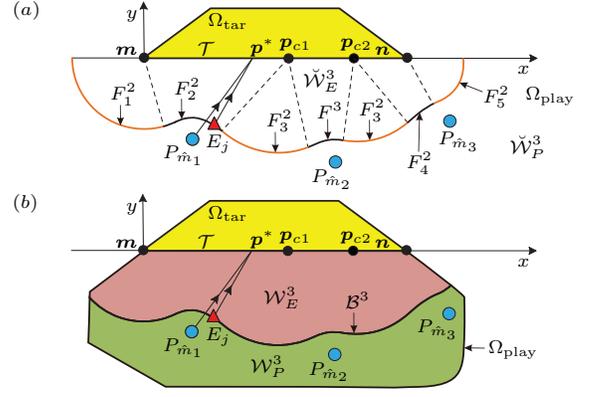}
\put(-143,15){$\scriptstyle{P_{\hat{m}_1}}$}
\put(-88,6){$\scriptstyle{P_{\hat{m}_2}}$}
\put(-47,22){$\scriptstyle{P_{\hat{m}_3}}$}
\put(-127,19){$\scriptstyle{E_j}$}
\put(-143,88){$\scriptstyle{P_{\hat{m}_1}}$}
\put(-87,78.9){$\scriptstyle{P_{\hat{m}_2}}$}
\put(-39,93){$\scriptstyle{P_{\hat{m}_3}}$}
\put(-127,92){$\scriptstyle{E_j}$}
\put(-110,7){$\scriptstyle{\mathcal{W}_P^3}$}
\put(-105,33){$\scriptstyle{\mathcal{W}_E^3}$}
\put(-74,31){$\scriptstyle{\mathcal{B}^3}$}
\put(-90,115){$\scriptstyle{\breve{\mathcal{W}}_E^3}$}
\put(-12,90){$\scriptstyle{\breve{\mathcal{W}}_P^3}$}
\put(-163,110){$\scriptstyle{F_1^2}$}
\put(-139,112){$\scriptstyle{F_2^2}$}
\put(-104,100){$\scriptstyle{F_3^2}$}
\put(-85,103){$\scriptstyle{F^3}$}
\put(-69,104){$\scriptstyle{F_3^2}$}
\put(-51,84){$\scriptstyle{F_4^2}$}
\put(-22,109){$\scriptstyle{F_5^2}$}
\put(-20,14){$\scriptstyle{\Omega_{\rm play}}$}
\put(-6,110){$\scriptstyle{\Omega_{\rm play}}$}
\put(-126,65){$\scriptstyle{\Omega_{\rm tar}}$}
\put(-126,138){$\scriptstyle{\Omega_{\rm tar}}$}
\put(-110,55){$\scriptstyle{\bm{p}^*}$}
\put(-110,128){$\scriptstyle{\bm{p}^*}$}
\put(-130,54){$\scriptstyle{\mathcal{T}}$}
\put(-130,127){$\scriptstyle{\mathcal{T}}$}
\put(-99,57){$\scriptstyle{\bm{p}_{c1}}$}
\put(-76,130){$\scriptstyle{\bm{p}_{c2}}$}
\put(-99,130){$\scriptstyle{\bm{p}_{c1}}$}
\put(-76,57){$\scriptstyle{\bm{p}_{c2}}$}
\put(-161,54){$\scriptstyle{\bm{m}}$}
\put(-161,127){$\scriptstyle{\bm{m}}$}
\put(-63,54){$\scriptstyle{\bm{n}}$}
\put(-63,127){$\scriptstyle{\bm{n}}$}
\put(-9,47){$\scriptstyle{x}$}
\put(-9,120){$\scriptstyle{x}$}
\put(-157,68){$\scriptstyle{y}$}
\put(-157,141){$\scriptstyle{y}$}
\put(-200,142){$\scriptstyle{(a)}$}
\put(-200,69){$\scriptstyle{(b)}$}
\caption{The barrier and two winning regions determined by three active pursuers $P_{\hat{m}_1},P_{\hat{m}_2}$ and $P_{\hat{m}_3}$ when playing with $E_j$, in which all three pursuers contribute to the construction of $\breve{\mathcal{B}}^3$. $(a)$ Without considering the boundary of $\Omega_{\rm play}$, the barrier $\breve{\mathcal{B}}^3$ consists of seven parts. Two winning regions $\breve{\mathcal{W}}_P^3$ and $\breve{\mathcal{W}}_E^3$ for each team's guaranteed winning are split by $\breve{\mathcal{B}}^3$. $(b)$ By considering the boundary of $\Omega_{\rm play}$, $\mathcal{B}^3$ is given by $\breve{\mathcal{B}}^3\cap\Omega_{\rm play}$, as the optimal trajectories for the players who work in the payoff function, are straight lines and $\Omega_{\rm play}$ is convex. Naturally, $\breve{\mathcal{W}}_P^3$ becomes the green region $\mathcal{W}_P^3$, and $\breve{\mathcal{W}}_E^3$ becomes the red region $\mathcal{W}_E^3$.\label{figT:5}}
\end{figure}

Next, a method is presented to construct the barrier $\breve{\mathcal{B}}^{n_k}$. First, by \autoref{barrier two}(b), construct the barrier $\breve{\mathcal{B}}^2$ determined by $P_{\hat{m}_1}$ and $P_{\hat{m}_2}$, as Fig. \ref{figT:4}(a) shows by replacing $P_{i1}$ and $P_{i2}$ with $P_{\hat{m}_1}$ and $P_{\hat{m}_2}$ respectively. Then, add the third pursuer $P_{\hat{m}_3}$ and compute the associated $\breve{\mathcal{B}}^3$ which is proved below to consist of seven parts, as Fig. \ref{figT:5}(a) shows.

According to the definition of the $x$-rank index set presented above, $x^0_{P_{\hat{m}_{1}}}<x^0_{P_{\hat{m}_{2}}}<x^0_{P_{\hat{m}_{3}}}$ holds. Take two points $\bm{p}_{c1}$ and $\bm{p}_{c2}$ in $\mathcal{T}$ such that\begin{equation}\label{pcv2v1definition}
\|\mathbf{x}_{P_{\hat{m}_i}}^0-\bm{p}_{ci}\|_2=\|\mathbf{x}_{P_{\hat{m}_{i+1}}}^0-\bm{p}_{ci}\|_2,i=1,2
\end{equation}
holds, as depicted in Fig. \ref{figT:5}(a), that is, $\bm{p}_{ci}$ is the point in $\mathcal{T}$ having the same distance to $P_{\hat{m}_i}$ and $P_{\hat{m}_{i+1}}$. Thus, from (\ref{pcv2v1definition}), $\bm{p}_{ci}=(x_{ci},0)$ can be computed as follows:\begin{equation}\label{pcv2v1}
x_{ci}=\frac{\|\mathbf{x}_{P_{\hat{m}_{i+1}}}^0\|_2^2-\|\mathbf{x}_{P_{\hat{m}_i}}^0\|_2^2}{2(x^0_{P_{\hat{m}_{i+1}}}-x^0_{P_{\hat{m}_i}})},i=1,2.
\end{equation}

Since $P_{\hat{m}_1}$ is an active pursuer, there must exist a point in $\mathcal{T}$ that $P_{\hat{m}_1}$ can reach before $P_{\hat{m}_2}$. Thus, $x_{c1}>0$, as Fig. \ref{figT:5}(a) shows. Similarly, $x_{c2}<l$ also holds. Take a point $\bm{p}=(x_p,0)$ in $\mathcal{T}$. It can be obtained that if $x_p\in[0,x_{c1})$, $P_{\hat{m}_1}$ can reach $\bm{p}$ before $P_{\hat{m}_2}$, and if $x_p\in(x_{c2},l]$, $P_{\hat{m}_3}$ can reach $\bm{p}$ before $P_{\hat{m}_2}$. Since $P_{\hat{m}_2}$ is an active pursuer, there must exist a point in $\mathcal{T}$ that $P_{\hat{m}_2}$ can reach before both $P_{\hat{m}_1}$ and $P_{\hat{m}_3}$. Thus, $x_{c1}<x_{c2}$ holds. In conclusion, $0<x_{c1}<x_{c2}<l$ is derived, as Fig. \ref{figT:5}(a) illustrates.

Assume $E_j$ can succeed to reach $\mathcal{T}$, and denote $E_j$'s OTP in $\mathcal{T}$ by $\bm{p}^*=(x_p^*,0)$ such that $J$ in (\ref{valueADR1v1}) only involving these three pursuers is maximized. Thus in this case,  (\ref{valueADR1v1}) becomes
\begin{equation}\begin{aligned}\label{valueADR3v1payoff}
J&=\min_{i=\hat{m}_1,\hat{m}_2,\hat{m}_3}\|\mathbf{x}_{P_i}(t_1)-\mathbf{x}_{E_j}(t_1)\|_2\\
V&=\min_{\mathbf{u}_{P_{\hat{m}_1}},\mathbf{u}_{P_{\hat{m}_2}}, \mathbf{u}_{P_{\hat{m}_3}}\in\mathcal{U}}\max_{\mathbf{u}_{E_j}\in\mathcal{U}}J
\end{aligned}\end{equation}
where $t_1$ is the first arrival time when $E_j$ reaches $\mathcal{T}$. Similarly, we further assume $\mathbf{x}_{E_j}^0\in\breve{\mathcal{B}}^3$. Therefore, under four players' optimal strategies from (\ref{valueADR3v1payoff}), $E_j$ will be captured exactly when reaching $\mathcal{T}$, namely, $V=0$.

 If $x_p^*\in[0,x_{c1})$, as Fig. \ref{figT:5}(a) shows, $P_{\hat{m}_1}$ can reach $\bm{p}^*$ before both $P_{\hat{m}_2}$ and $P_{\hat{m}_3}$. Thus, this part of $\breve{\mathcal{B}}^{3}$ is determined by $P_{\hat{m}_1}$ alone, namely, $J$ only depends on $P_{\hat{m}_1}$. By following the proof of \autoref{barrier two}, the part related to $P_{\hat{m}_1}$ alone is given by $F^2_1(\mathbf{x}^0_{P_{\hat{m}_1}},\mathbf{x}^0_{P_{\hat{m}_2}})\cup F_2^2(\mathbf{x}^0_{P_{\hat{m}_1}},\mathbf{x}^0_{P_{\hat{m}_{2}}})$, respectively for $x_p^*=0$ and $x_p^*\in(0,x_{c1})$, which consists of the leftmost orange and black curves in Fig. \ref{figT:5}(a).
 
 If $x_p^*\in(x_{c1},x_{c2})$, as Fig. \ref{figT:5}(a) shows, $P_{\hat{m}_2}$ can reach $\bm{p}^*$ before both $P_{\hat{m}_1}$ and $P_{\hat{m}_3}$, implying that this part of $\breve{\mathcal{B}}^{3}$ is determined by $P_{\hat{m}_2}$ alone, namely, $J$ only depends on $P_{\hat{m}_2}$. Thus, it follows from the proof of \autoref{barrierone} that (\ref{OTPunequ11}) and (\ref{barrier1v1explicit}) hold for $P_{\hat{m}_2}$ and $E_j$, that is, $x_p^*=\frac{x_{E_j}^0-\alpha^2x_{P_{\hat{m}_2}}^0}{1-\alpha^2}$ and
\begin{equation}\begin{aligned}\label{barrier3v1explicitDim2}
(x_{E_j}^0-x_{P_{\hat{m}_2}}^0)^2+(1-1/\alpha^2)(y_{E_j}^0)^2+(1-\alpha^2)(y_{P_{\hat{m}_2}}^0)^2 =0
\end{aligned}\end{equation} 
which characterizes the relationship between the initial positions of $P_{\hat{m}_2}$ and $E_j$ when $\mathbf{x}_{E_j}^0\in\breve{\mathcal{B}}^3$ and $x_p^*\in(x_{c1},x_{c2})$. Also note that from $x_p^*\in(x_{c1},x_{c2})$, we have 
 \begin{equation}\begin{aligned}\label{barrier3v1intervalb33}
x_{E_j}^0\in\big((1-\alpha^2)x_{c1}+\alpha^2x_{P_{\hat{m}_2}}^0,(1-\alpha^2)x_{c2}+\alpha^2x_{P_{\hat{m}_2}}^0\big).
\end{aligned}\end{equation}
Therefore, by replacing $\mathbf{x}_{E_j}^0$ with a general variable $\mathbf{z}\in\mathbb{R}^2$, it follows from (\ref{barrier3v1explicitDim2}), (\ref{barrier3v1intervalb33}) and $y_{E_j}^0<0$ that this part of $\breve{\mathcal{B}}^{3}$, which is the middle black curve in Fig. \ref{figT:5}(a), is given by
\begin{equation}\begin{aligned}\label{barrier2v1B22}
&\big\{\mathbf{z}=(x,y)\in\mathbb{R}^2|(x-x_{P_{\hat{m}_2}}^0)^2+(1-1/\alpha^2)y^2\\
&\qquad +(1-\alpha^2)(y_{P_{\hat{m}_2}}^0)^2 =0,x\in(k_7,k_8),y<0\big\}
\end{aligned}\end{equation}
where $k_7=(1-\alpha^2)x_{c1}+\alpha^2x_{P_{\hat{m}_2}}^0$ and $k_8=(1-\alpha^2)x_{c2}+\alpha^2x_{P_{\hat{m}_2}}^0$. From (\ref{basecurvethree}), it can be noted that (\ref{barrier2v1B22}) can be expressed by a base curve, that is, $F^3(\mathbf{x}^0_{P_{\hat{m}_1}},\mathbf{x}^0_{P_{\hat{m}_{2}}},\mathbf{x}^0_{P_{\hat{m}_3}})$.

If $x_p^*\in(x_{c2},l]$, as Fig. \ref{figT:5}(a) shows, this part of $\breve{\mathcal{B}}^3$ is only related to $P_{\hat{m}_3}$. Based on the similar analysis for the case $x_p^*\in[0,x_{c1})$, this part, which consists of the rightmost black and orange curves in Fig. \ref{figT:5}(a), is given by $F_4^2(\mathbf{x}^0_{P_{\hat{m}_2}},\mathbf{x}^0_{P_{\hat{m}_3}})\cup F_5^2(\mathbf{x}^0_{P_{\hat{m}_2}},\mathbf{x}^0_{P_{\hat{m}_3}})$. 

Now, we consider the remainder $x_p^*=x_{c1}$ or $x_p^*=x_{c2}$. For simplicity, we only consider $x_p^*=x_{c1}$, and similar analysis can be conducted for $x_p^*=x_{c2}$. Note that $x_p^*=x_{c1}$, i.e., $\bm{p}^*=\bm{p}_{c1}$, and (\ref{distanceunequ}) holds for $P_{\hat{m}_1}$ and $E_j$. Similar to the analysis for $x_p^*=x_c$ in the proof of \autoref{barrier two}, this part of $\breve{\mathcal{B}}^3$, which is the second orange curve from the left in Fig. \ref{figT:5}(a), is given by a base curve in (\ref{basecurves22}), that is, $F^2_3(\mathbf{x}^0_{P_{\hat{m}_1}},\mathbf{x}^0_{P_{\hat{m}_{2}}})$.

Thus, by connecting all these parts of $\breve{\mathcal{B}}^{3}$, we can obtain $\breve{\mathcal{B}}^{3}$ given by (\ref{barriergeneralnv1}) with $n_k=3$, shown in Fig. \ref{figT:5}(a).

Similarly, add the fourth pursuer $P_{\hat{m}_4}$, and $\bm{p}_{c3}=(x_{c3},0)$ is given by (\ref{pcv2v1}) with $i=3$. Then, $0<x_{c1}<x_{c2}<x_{c3}<l$ can be obtained. In the same way, $\breve{\mathcal{B}}^{4}$ is given by (\ref{barriergeneralnv1}) with $n_k=4$. Therefore, by adding the remaining pursuers one by one along the $x$-rank index set $\mathcal{\hat{I}}_k$, $\breve{\mathcal{B}}^{n_k}$ is obtained as (\ref{barriergeneralnv1}) shows. 

Then, consider the effect of the boundary of $\Omega_{\rm play}$. Since the optimal trajectories for the players who work in the payoff function $J$, are straight lines and $\Omega_{\rm play}$ is convex, it can be obtained that $\mathcal{B}^{n_k}=\breve{\mathcal{B}}^{n_k}\cap\Omega_{\rm play}$, as Fig. \ref{figT:5}(b) shows when $n_k=3$. Naturally, the PWR $\mathcal{W}_P^3$ is the green region, and the EWR $\mathcal{W}_E^3$ is the red region. Thus we finish the proof.\qed

Note that in general pursuit coalitions, there may exist some inactive pursuers who as asserted below, have no effect on the barrier construction but instead hinder our analysis. In light of this, if all active pursuers can be extracted out from any given pursuit coalition, namely, its largest full-active pursuit subcoalition, then the barrier for this general pursuit coalition can be obtained from \autoref{barrierfull}. 

First, the following lemma formally establishes a connection of the barrier between a general pursuit coalition and its largest full-active pursuit subcoalition.

\begin{lema}[Barrier Equivalence]\label{barriereqal}\rm
For any pursuit coalition $k$, $\mathcal{B}^{n_k}(\mathcal{X}_k^0)=\mathcal{B}^{\bar{n}_k}(\mathcal{\bar{X}}_k^0)$ holds, where $\mathcal{\bar{X}}_k^0$ with $\bar{n}_k$ pursuers denotes the initial position of the unique largest full-active pursuit subcoalition of $\mathcal{X}_k^0$.
\end{lema}

\emph{Proof:} Obviously, we only need to focus on the case that $\mathcal{X}_k^0\setminus\mathcal{\bar{X}}^0_k$ is nonempty. Thus, $n_k>\bar{n}_k$. For simplicity, we prove $\mathcal{W}_E^{n_k}=\mathcal{W}_E^{\bar{n}_k}$. By definition, $\mathcal{W}_E^{n_k}\subseteq\mathcal{W}_E^{\bar{n}_k}$ holds, as reducing the pursuer will result in the same or expansion of the EWR. Thus we only need to verify $\mathcal{W}_E^{\bar{n}_k}\subseteq\mathcal{W}_E^{n_k}$, equivalently, $\breve{\mathcal{W}}_E^{\bar{n}_k}\subseteq\breve{\mathcal{W}}_E^{n_k}$.

Suppose $\bm{p}\in\breve{\mathcal{W}}_E^{\bar{n}_k}$. Then, there must exist a point $\bm{p}_1\in\mathcal{T}$ such that\begin{equation}\label{barrierequa11}
\|\bm{p}-\bm{p}_1\|_2<\alpha \|\mathbf{x}_{P_i}^0-\bm{p}_1\|_2
\end{equation}
holds for all $\mathbf{x}_{P_i}^0\in\mathcal{\bar{X}}_k^0$. Assume that there exists a pursuer $P_j$ satisfying $\mathbf{x}_{P_j}^0\in\mathcal{X}_k^0\setminus\mathcal{\bar{X}}^0_k$ and \begin{equation}\label{barrierequa22}
\alpha \|\mathbf{x}_{P_j}^0-\bm{p}_1\|_2\leq\|\bm{p}-\bm{p}_1\|_2
\end{equation}
Then, by (\ref{barrierequa11}) and (\ref{barrierequa22}), it can be stated that\begin{equation}\label{barrierequa33}
\|\mathbf{x}_{P_j}^0-\bm{p}_1\|_2<\|\mathbf{x}_{P_i}^0-\bm{p}_1\|_2
\end{equation}
holds for all $\mathbf{x}_{P_i}^0\in\mathcal{\bar{X}}_k^0$. Thus, there exists a point, i.e, $\bm{p}_1$, in $\mathcal{T}$ that $P_j$ can reach before all pursuers in $\mathcal{\bar{X}}_k^0$, which contradicts with the fact that $\mathcal{\bar{X}}_k^0$ is the largest full-active pursuit subcoalition of $\mathcal{X}_k^0$. Therefore, we can conclude that (\ref{barrierequa22}) does not hold, that is,
\begin{equation}\label{barrierequa44}
\|\bm{p}-\bm{p}_1\|_2<\alpha \|\mathbf{x}_{P_j}^0-\bm{p}_1\|_2
\end{equation}
holds for all $\mathbf{x}_{P_j}^0\in\mathcal{X}_k^0\setminus\mathcal{\bar{X}}^0_k$. Combining (\ref{barrierequa11}) with (\ref{barrierequa44})  shows that there exists a point, i.e., $\bm{p}_1$, in $\mathcal{T}$ that $\bm{p}$ can reach before all pursuers in $\mathcal{X}_k^0$ with a speed ratio $\alpha$. Thus, $\bm{p}\in\breve{\mathcal{W}}_E^{n_k}$, implying that $\breve{\mathcal{W}}_E^{\bar{n}_k}\subseteq\breve{\mathcal{W}}_E^{n_k}$.

The existence of $\mathcal{\bar{X}}_k^0$ is straightforward, as $\mathcal{\bar{X}}_k^0$ is a subset of $\mathcal{X}_k^0$. As for the uniqueness, assume that $\mathcal{\bar{X}}_k^0$ with $\bar{n}_k$ pursuers and $\mathcal{\bar{Y}}_k^0$ with $\bar{m}_k$ pursuers are two distinct largest full-active pursuit subcoalitions of $\mathcal{X}_k^0$.

Take a pursuer $P_i$ satisfying $\mathbf{x}_{P_i}^0\in\mathcal{\bar{X}}_k^0$ and $\mathbf{x}_{P_i}^0\notin\mathcal{\bar{Y}}_k^0$. Then, it follows from \autoref{full active} that $\mathbf{x}_{P_i}^0\notin\mathcal{\bar{Y}}_k^0$ implies that there is no point in $\mathcal{T}$ that $P_i$ can reach before all pursuers in $\mathcal{\bar{Y}}_k^0$. In other words, there is no point in $\mathcal{T}$ that $P_i$ can reach before all other pursuers in $\mathcal{X}_k^0$, which contradicts with the fact that $\mathbf{x}_{P_i}^0$ is an active pursuer by noting that $\mathbf{x}_{P_i}^0\in\mathcal{\bar{X}}_k^0$. Thus, $\mathcal{\bar{X}}_k^0=\mathcal{\bar{Y}}_k^0$ and the uniqueness is proved.\qed

Denote by $\mathcal{R}_D$ the set of points in $\mathbb{R}^2$ that $P_{m_i}$ can reach before $P_{m_j}$, which is mathematically formulated as follows:
\begin{equation}\begin{aligned}\label{pursuerdominancere}
\mathcal{R}_D\big(\mathbf{x}_{P_{m_i}}^0,\mathbf{x}_{P_{m_j}}^0\big)=\big\{\mathbf{z}\in\mathbb{R}^2|\|\mathbf{z}&-\mathbf{x}_{P_{m_i}}^0\|_2\\
&<\|\mathbf{z}-\mathbf{x}_{P_{m_j}}^0\|_2\big\}.
\end{aligned}\end{equation}
It can be noted that $\mathcal{R}_D$ is a half plane.

Next, for a pursuit coalition $k$, an algorithm is presented to find its unique largest full-active pursuit subcoalition, namely, $\mathcal{\bar{X}}_k^0$. To avoid redundant illustration, the proof of this algorithm is left in the next theorem.

\begin{algm}\label{algm1}\rm
(Algorithm for Finding the Largest Full-Active Pursuit Subcoalition in $k$).\\
1: \textbf{Input:} $\mathcal{X}_k^0\in\Omega_{\rm play}^{n_k}\cup\mathcal{T}^{n_k},\bar{\mathcal{X}}_k^0\gets\emptyset,\bar{n}_k\gets0$.\\
2: \textbf{for} $i\in\{1,...,n_k\}$ \textbf{do}\\
3: \hspace{0.2 in}$\mathcal{T}_1=\mathcal{T}$\\
4: \hspace{0.2 in}\textbf{for} $j\in\{1,...,n_k\},j\neq i$ \textbf{do}\\
5: \hspace{0.4 in}$\mathcal{T}_1=\mathcal{R}_D(\mathbf{x}_{P_{m_i}}^0,\mathbf{x}_{P_{m_j}}^0)\cap\mathcal{T}_1$ \textbf{end for}\\
6: \hspace{0.2 in}\textbf{if} $\mathcal{T}_1\neq\emptyset$ \textbf{then}\\
7: \hspace{0.4 in}$\bar{\mathcal{X}}_k^0\gets\bar{\mathcal{X}}_k^0\cup \mathbf{x}^0_{P_{m_i}},\bar{n}_k\gets\bar{n}_k+1$ \textbf{end if} \textbf{end for}\\
8: \textbf{Output:} $\bar{\mathcal{X}}_k^0,\bar{n}_k$.
\end{algm}

Now it suffices to present the main result of this section.

\begin{thom}[Barrier for General Pursuit Coalition]\label{barriergeneral}\rm
Consider the system (1) satisfying Assumptions \ref{isolateinitial}, \ref{constrainedinitial} and \ref{speedratio}. For a pursuit coalition $k$, its largest full-active pursuit subcoalition $\mathcal{\bar{X}}_k^0$ with $\bar{n}_k$ pursuers can be found by \autoref{algm1}. Then, $\mathcal{B}^{n_k}(\mathcal{X}_k^0)=\mathcal{B}^{\bar{n}_k}(\mathcal{\bar{X}}_k^0)$, where $\mathcal{B}^{\bar{n}_k}(\mathcal{\bar{X}}_k^0)$ can be computed by \autoref{barrierfull}.
\end{thom}
 
\emph{Proof:} According to \autoref{barriereqal} and \autoref{barrierfull}, if the validity of \autoref{algm1} can be verified, this theorem holds naturally. 

From line 3 to line 5 in \autoref{algm1}, the set of points in $\mathcal{T}$ that $P_{m_i}$ can reach before all other pursuers in $k$, is computed and denoted by $\mathcal{T}_1$, where $\mathcal{R}_D(\mathbf{x}_{P_{m_i}}^0,\mathbf{x}_{P_{m_j}}^0)$ denotes the set of points in $\mathbb{R}^2$ that $P_{m_i}$ can reach before $P_{m_j}$ as (\ref{pursuerdominancere}) shows. Thus, if $\mathcal{T}_1\neq\emptyset$, that is, there exists at least one point that $P_i$ can reach before all other pursuers in $k$, then $P_i$ must be an active pursuer in $k$. Conversely, if $\mathcal{T}_1=\emptyset$, namely, there is no point in $\mathcal{T}$ that $P_i$ can reach before all other pursuers in $k$, then $P_i$ must be an inactive pursuer in $k$.\qed

\subsection{Extensions to Relaxed Initial Deployment}\label{sectionextension}
In this section, the barrier construction method will be extended to the scenarios having more realistic initial deployment as \autoref{relaxedinitial} shows. In this initial deployment, pursuers are allowed to initially lie in $\Omega_{\rm tar}$. For example, there may exist pursuers who are patrolling in $\Omega_{\rm tar}$ exactly when evaders are detected. As discussed below, by introducing a specific virtual pursuer, a road between this case and the known results is established. Or, more concretely, it is demonstrated that this virtual pursuer plays the same role with its original pursuer in barrier construction.

\begin{defi}[Virtual Pursuer]\label{virtual}\rm
For every pursuer $P_i$ whose initial condition satisfies $\mathbf{x}_{P_i}^0=(x_{P_i}^0,y_{P_i}^0)\in\Omega_{\rm tar}$, introduce a virtual pursuer $\tilde{P}_i$ with the initial position $\tilde{\mathbf{x}}_{P_i}^0=(\tilde{x}_{P_i}^0,\tilde{y}_{P_i}^0)$ such that $\tilde{x}_{P_i}^0=x_{P_i}^0$ and $\tilde{y}_{P_i}^0=-y_{P_i}^0$.
\end{defi}

\begin{figure}\centering
\graphicspath{{Figures/}}
\includegraphics[width=77mm,height=25mm]{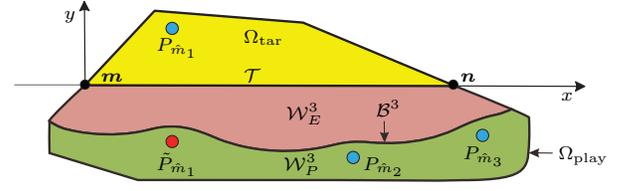}
\put(-163,51){$\scriptstyle{P_{\hat{m}_1}}$}
\put(-163,6){$\scriptstyle{\tilde{P}_{\hat{m}_1}}$}
\put(-85,6){$\scriptstyle{P_{\hat{m}_2}}$}
\put(-47,10){$\scriptstyle{P_{\hat{m}_3}}$}
\put(-115,6){$\scriptstyle{\mathcal{W}_P^3}$}
\put(-114,25){$\scriptstyle{\mathcal{W}_E^3}$}
\put(-80,26){$\scriptstyle{\mathcal{B}^3}$}
\put(-11,10){$\scriptstyle{\Omega_{\rm play}}$}
\put(-130,54){$\scriptstyle{\Omega_{\rm tar}}$}
\put(-130,39.5){$\scriptstyle{\mathcal{T}}$}
\put(-184,39.5){$\scriptstyle{\bm{m}}$}
\put(-48,39.5){$\scriptstyle{\bm{n}}$}
\put(-10,32){$\scriptstyle{x}$}
\put(-198,63){$\scriptstyle{y}$}
\caption{The barrier and winning regions for relaxed initial deployment which allows the pursuer to start the game from $\Omega_{\rm tar}$. The red circle $\tilde{P}_{\hat{m}_1}$ is the virtual pursuer of $P_{\hat{m}_1}$, and they are symmetric with respect to $\mathcal{T}$. It is proved that the barrier $\mathcal{B}^3$ determined by $P_{\hat{m}_1},P_{\hat{m}_2}$ and $P_{\hat{m}_3}$ is the same as the barrier determined by $\tilde{P}_{\hat{m}_1},P_{\hat{m}_2}$ and $P_{\hat{m}_3}$. This feature is called mirror property.\label{figBT:1}}
\end{figure}

\begin{rek}\label{rek34}\rm
It can be easily observed that the virtual pursuer and its original pursuer are symmetric with respect to $\mathcal{T}$, and a three-pursuer pursuit coalition example is presented in Fig. \ref{figBT:1} where $P_{\hat{m}_1}$ lies in $\Omega_{\rm tar}$ and its virtual pursuer is the red circle $\tilde{P}_{\hat{m}_1}$. Since $\Omega_{\rm play}$ and $\Omega_{\rm tar}$ allow different shapes, the virtual pursuer may lie out of $\Omega$. Recalling the barrier construction process under Assumptions \ref{isolateinitial} and \ref{constrainedinitial}, it can be claimed that this case can also be unified by constructing $\breve{\mathcal{B}}^{n_k}$ first and then computing $\breve{\mathcal{B}}^{n_k}\cap\Omega_{\rm play}$, as the optimal trajectories for the players are straight lines and $\Omega$ is convex.
\end{rek}

\begin{rek}\label{rek35}\rm
Note that the virtual pursuer generated by one pursuer may coincide with another pursuer, which will lead to a disagreement with \autoref{isolateinitial}. Thus, assume that every virtual pursuer does not coincide with all other pursuers. However, the coincidence of the virtual pursuer and its own original pursuer is allowed, namely, when $\mathbf{x}_{P_i}^0\in\mathcal{T}$.
\end{rek}

Next, an important property of virtual pursuers in terms of the barrier construction will be stated, which is a key step to simplify the problems.

\begin{lema}[Mirror Property]\label{mirror}\rm
For a pursuit coalition $k$, if $\mathbf{x}_{P_i}^0\in\Omega_{\rm tar}$ $(i\in\mathcal{I}_k)$, let $\mathcal{X}_k^0(-i)$ denote the remaining pursuers in $k$ when $\mathbf{x}_{P_i}^0$ is eliminated, and $\tilde{\mathbf{x}}_{P_i}^0$ denote the virtual pursuer of $\mathbf{x}_{P_i}^0$. Then, $\mathcal{B}^{n_k}(\mathcal{X}^0_k)=\mathcal{B}^{n_k}(\mathcal{X}^0_k(-i)\cup\tilde{\mathbf{x}}_{P_i}^0)$.
\end{lema}

\emph{Proof:} Suppose $\bm{p}\in\mathcal{W}_E^{n_k}(\mathcal{X}^0_k(-i)\cup\tilde{\mathbf{x}}_{P_i}^0)$, and then there must exist a point $\bm{p}_1\in\mathcal{T}$ such that 
\begin{equation}\label{mirriorequ1}
\|\bm{p}-\bm{p}_1\|_2<\alpha \|\mathbf{x}_{P_j}^0-\bm{p}_1\|_2 
\end{equation}
holds for all $\mathbf{x}_{P_j}^0\in\mathcal{X}^0_k(-i)\cup\tilde{\mathbf{x}}_{P_i}^0$. According to \autoref{virtual}, we have
\begin{equation}\label{mirriorequ2}
\|\mathbf{x}_{P_i}^0-\bm{p}_1\|_2=\|\tilde{\mathbf{x}}_{P_i}^0-\bm{p}_1\|_2.
\end{equation}
Naturally, (\ref{mirriorequ1}) holds for all $\mathbf{x}_{P_j}^0\in\mathcal{X}^0_k$, which implies that $\bm{p}\in\mathcal{W}_E^{n_k}(\mathcal{X}^0_k)$. Therefore, $\mathcal{W}_E^{n_k}(\mathcal{X}^0_k(-i)\cup\tilde{\mathbf{x}}_{P_i}^0)\subset\mathcal{W}_E^{n_k}(\mathcal{X}^0_k)$ is obtained.

On the other side, suppose $\bm{p}\in\mathcal{W}_E^{n_k}(\mathcal{X}^0_k)$ and in the similar way, $\bm{p}\in\mathcal{W}_E^{n_k}(\mathcal{X}^0_k(-i)\cup\tilde{\mathbf{x}}_{P_i}^0)$ can be derived, implying $\mathcal{W}_E^{n_k}(\mathcal{X}^0_k)\subset\mathcal{W}_E^{n_k}(\mathcal{X}^0_k(-i)\cup\tilde{\mathbf{x}}_{P_i}^0)$. Thus,  $\mathcal{W}_E^{n_k}(\mathcal{X}^0_k)=\mathcal{W}_E^{n_k}(\mathcal{X}^0_k(-i)\cup\tilde{\mathbf{x}}_{P_i}^0)$, which finishes the proof by noting that $\mathcal{B}^{n_k}$ is the separating curve between $\mathcal{W}_P^{n_k}$ and $\mathcal{W}_E^{n_k}$. A three pursuer case is presented in Fig. \ref{figBT:1}, where $P_{\hat{m}_1}$ lies in $\Omega_{\rm tar}$.\qed

The main result in this section is now given, which provides an efficient way to construct the barrier under Assumptions \ref{isolateinitial} and \ref{relaxedinitial} by projecting this problem into the field of the former section. 

\begin{corol}[Barrier for Relaxed Initial Deployment]\label{barrier relaxed}\rm
Consider the system (1) and suppose that Assumptions \ref{isolateinitial}, \ref{relaxedinitial} and \ref{speedratio} hold. For a pursuit coalition $k$, let $\mathcal{X}_{k,1}^0$ and $\mathcal{X}_{k,2}^0$ (maybe empty) denote the initial sets of its pursuers who lie in $\Omega_{\rm play}$ and $\Omega_{\rm tar}$ respectively. By one-to-one mapping of $\mathcal{X}_{k,2}^0$, a set of virtual pursuers can be obtained and denote it by $\mathcal{\tilde{X}}_{k,2}^0$. Then, $\mathcal{B}^{n_k}(\mathcal{X}^0_k)=\mathcal{B}^{n_k}(\mathcal{X}_{k,1}^0\cup\mathcal{\tilde{X}}_{k,2}^0)$ holds, where $\mathcal{B}^{n_k}(\mathcal{X}_{k,1}^0\cup\mathcal{\tilde{X}}_{k,2}^0)$ can be computed by \autoref{barriergeneral}.
\end{corol} 

\emph{Proof:} Since all players in $\mathcal{X}_{k,1}^0\cup\mathcal{\tilde{X}}_{k,2}^0$ lie in $\Omega_{\rm play}\cup\mathcal{T}$, this corollary is straightforward by considering \autoref{mirror} and \autoref{barriergeneral}.\qed

\section{Pursuit Task Assignment}\label{task}
In the next, by matching pursuit coalitions with evaders, an optimal task assignment scheme for the pursuit team to guarantee the most evaders intercepted, will be investigated. Intuitively, for every evader, we want to designate a pursuit coalition which can make sure of capturing it. In view of the characteristics of winning regions, the selection of an adequate pursuit coalition can be realized by checking if this evader lies in the capturable region of this pursuit coalition, namely, $\mathcal{W}_P^{n_k}$. In this way, rich prior information about which evaders can be captured by a specified pursuit coalition, is collected. Then, pursuit coalitions can be matched with evaders one by one such that the most evaders are captured. Finally, this maximum matching is formulated as a simplified 0-1 integer programming problem instead of solving a constrained bipartite matching problem \cite{lawler1976combinatorial}. 

Let $G=(P,E,\mathcal{E})$ denote an undirected bipartite graph, consisting of two independent sets $P,E$ of nodes, and a set $\mathcal{E}$ of unordered pairs of nodes called edges each of which connects a node in $P$ to one in $E$. In this case, we take $P=\{1,2,...,2^{N_p}-1\}$, representing all pursuit coalitions, and $E=\{1,2,..,N_e\}$ on behalf of all evaders. An edge from node $i$ in $P$ to node $j$ in $E$ is denoted by $e_{ij}$, referring to using pursuit coalition $i$ to intercept evader $E_j$. Define $e_{ij}\in \mathcal{E}$ if pursuit coalition $i$ can guarantee the capture of evader $E_j$ in $\Omega_{\rm play}$ or $\mathcal{T}$; otherwise, $e_{ij}\notin \mathcal{E}$. Thus, by computing the barrier for every pursuit coalition, all edges contained in $\mathcal{E}$ can be found.

Traditionally, a subset $M\subseteq \mathcal{E}$ is said to be a matching if no two edges in $M$ are incident to the same node, and our goal is to find a matching containing a maximum number of edges. However, note that every pursuer can appear in at most one pursuit coalition when the interception scheme is executed. Thus, the pursuit coalitions containing at least one same pursuer cannot coexist, which results in a constrained bipartite matching problem. To solve it, this problem is transformed into the framework of 0-1 integer programming. 

Before proceeding with this transformation, the following lemma is presented, which will dramatically decrease the complexity of the 0-1 integer programming proposed later.

\begin{lema}[Degeneration of Pursuit Coalition]\label{degene}\rm
For any pursuit coalition $k$ with $n_k\ge3$, if there exists an evader $E_j$ such that $\mathbf{x}_{E_j}^0\in\mathcal{W}_P^{n_k}$, then there must exist a pursuit subcoalition $k_1$ of $k$ satisfying $n_{k_1}=2$ and $\mathbf{x}_{E_j}^0\in\mathcal{W}_P^{n_{k_1}}$.
\end{lema}

\emph{Proof:} Note that \autoref{barrierfull} manifests a special feature of the barrier that its every part is only associated with at most two pursuers as (\ref{barriergeneralnv1}) shows. Even for the part $F^3$, one can split it into two parts both of which are only associated with two pursuers. Thus, if a pursuit coalition $k$ can guarantee to capture an evader, then there must exist an its two-pursuer subcoalition which can also guarantee the capture of this evader.\qed

\begin{rek}\label{rek41}\rm
From \autoref{degene}, it can be stated that any maximum matching can be reduced to a simpler version in which every pursuit coalition contains at most two pursuers. Therefore, seeking for a maximum matching in the class of this simple version suffices to obtain a matching that is global optimal in the sense of the most evaders intercepted. 
\end{rek}

\begin{rek}\label{rekadd41}\rm
Although it follows from \autoref{degene} and \autoref{rek41} that the pursuit coalitions with more than two pursuers are not necessary in the maximum matching, the analysis for general pursuit coalitions in Section \ref{sectionmultiplev1} can provide instructions to determine the capturable and uncapturable regions of multiple pursuers with any numbers when pursuing one evader. Additionally, this result can also help to deploy the pursuers such that their capturable region is desirable, as sometimes the pursuit team should position its members before an evader occurs. More importantly, it is \autoref{barrierfull} that reveals the degeneration property of the pursuit coalition described in \autoref{degene}.
\end{rek}

Thus, in the following discussion, attention will be focused on this specific class of pursuit coalitions with at most two pursuers, and due to its crucial role in extracting out an optimal task assignment, it is stated formally as follows.

\begin{defi}[Execution Pursuit Coalition]\label{execution}\rm
A pursuit coalition $k$ is called an execution pursuit coalition, if $n_k=1$ or $n_k=2$.
\end{defi}

According to the barrier construction in Section \ref{pursuitcoalition}, the following prior information vector can be acquired. For notational convenience, define $N_v=N_eN_p(N_p+1)/2$.

\begin{defi}[Prior Information Vector]\label{prior}\rm
For $P_i$, define $\bm{r}_i^1=[r_i^1(1),...,r_i^1(N_e)]\in\mathbb{R}^{N_e}$, where for $j=1,...,N_e$, set $r_i^1(j)=0$ if $\mathbf{x}^0_{E_j}\in\mathcal{W}^1_E(\mathbf{x}^0_{P_i})$, that is, $P_i$ cannot guarantee to capture $E_j$ prior to its arrival in $\mathcal{T}$; otherwise, set $r_i^1(j)=1$. Similarly, for $P_{i1}$ and $P_{i2}$, define $\bm{r}_{i1,i2}^2=[r_{i1,2}^2(1),...,r_{i1,i2}^2(N_e)]\in\mathbb{R}^{N_e}$, where for $j=1,...,N_e$, set $r_{i1,i2}^2(j)=0$ if $\mathbf{x}^0_{E_j}\in\mathcal{W}^2_E(\mathbf{x}^0_{P_{i1}}\cup\mathbf{x}^0_{P_{i2}})$; otherwise, set $r_{i1,i2}^2(j)=1$. Then we call this vector $\bm{r}=[\bm{r}_1^1,...,\bm{r}_{N_p}^1,\bm{r}_{1,2}^2,...\bm{r}_{1,N_p}^2,\bm{r}^2_{2,3},...,\bm{r}_{N_p-1,N_p}^2]\in\mathbb{R}^{N_v}$ as prior information vector. 
\end{defi}

\begin{rek}\label{rek42}\rm
Note that the prior information vector $\bm{r}$ contains all information by which for any given evader, one can judge whether an execution pursuit coalition can guarantee to capture it. This vector is the only input for the maximum matching.
\end{rek}

Let $\bm{s}_i^1=[s_i^1(1),...,s_i^1(N_e)]\in\mathbb{R}^{N_e}$ denote the strategy vector of $P_i$. Its elements are either 1 or 0, where $s^1_i(j)=1$ indicates the assignment of $P_i$ to intercept $E_j$, and $s^1_i(j)=0$ means no assignment. Clearly, $\sum_{j=1}^{N_e}s^1_i(j)\leq1$ must be satisfied, namely, $P_i$ at a time can pursue at most one evader. Let $\bm{s}_{i1,i2}^2=[s_{i1,i2}^2(1),...,s_{i1,i2}^2(N_e)]\in\mathbb{R}^{N_e}$ denote the strategy vector of the pursuit pair $\{P_{i1},P_{i2}\}$. Its elements are either 1 or 0. Specifically, $s_{i1,i2}^2(j)=1$ indicates that the pursuit pair $\{P_{i1},P_{i2}\}$ cooperates to intercept $E_j$, and $s_{i1,i2}^2(j)=0$ means no assignment. Obviously, $\sum_{j=1}^{N_e}s_{i1i2}^2(j)\leq1$.

Denote the joint strategy vector of all one-pursuer pursuit coalitions by $\bm{s}^1=[\bm{s}_1^{1},...,\bm{s}_{N_p}^{1}]\in\mathbb{R}^{N_pN_e}$, and denote $\bm{s}^2=[\bm{s}_{1,2}^{2},...,\bm{s}^2_{1,N_p},\bm{s}_{2,3}^{2},...,\bm{s}_{N_p-1,N_p}^{2}]\in\mathbb{R}^{N_eN_p(N_p-1)/2}$  as the joint strategy vector of all two-pursuer pursuit coalitions. Thus $\bm{z}=[\bm{s}^1,\bm{s}^2]^\mathsf{T}\in\mathbb{R}^{N_v\times 1}$ denotes all execution pursuit coalitions' strategy vector.

Let ${\rm ones}(m,n)$ denote the $m\times n$ matrix each element of which is 1, ${\rm zeros}(m,n)$ denote the $m\times n$ zero matrix, $I_n$ denote the identity matrix of size $n$, $\otimes$ denote the Kronecker product, and $\mathcal{X}^0_P,\mathcal{X}^0_E$ denote the initial positions of all pursuers and all evaders, respectively.

Now, the main result of this section is presented below, which gives a maximum matching solution by solving a 0-1 integer programming.

\begin{thom}[Maximum Matching]\label{maximum}\rm
Consider the system (1) and suppose that Assumptions \ref{isolateinitial}, \ref{relaxedinitial} and \ref{speedratio} hold. Given $\mathcal{X}^0_P$ and $\mathcal{X}^0_E$, the number $q$ of the evaders which the pursuit team can guarantee to prevent from reaching the target region $\Omega_{\rm tar}$ is given by
\begin{equation}\begin{aligned}\label{maximizeequal11}
q&=\text{maximize }\bm{c}^\mathsf{T}\bm{z}\qquad\qquad\qquad\qquad\qquad\qquad\\
&\quad\ \text{s. t. }A_1\bm{z}\leq \bm{b}_1,A_2\bm{z}\leq \bm{b}_2,A_3\bm{z}\leq \bm{b}_3\\
&\quad\ \bm{z}=[\bm{s}^1,\bm{s}^2]^\mathsf{T}=[z(1),...,z(N_v)]^\mathsf{T},z(i)=0,1
\end{aligned}\end{equation}
and the maximum matching is given by $\bm{z}^*={\rm argmax}_{\bm{z}}(\bm{c}^\mathsf{T}\bm{z})$. The parameter matrixes and vectors are defined as follows:
\begin{equation}\begin{aligned}
&\bm{c}={\rm ones}(N_v,1),\bm{b}_1=\bm{r}^\mathsf{T},A_1=I_{N_v},\bm{b}_2={\rm ones}(N_e,1),\\
&A_2={\rm ones}(1,N_v/N_e)\otimes I_{N_e},\bm{b}_3={\rm ones}(N_p,1)
\end{aligned}\end{equation}
and $A_3$ is computed by \autoref{matrix} presented below.
\end{thom}

\emph{Proof:} Note that the objective function is the number of evaders which are assigned execution pursuit coalitions, that is, the number of the evaders to be captured. The first inequality constraint in (\ref{maximizeequal11}) represents that the prior information must be satisfied. In other words, assign an adequate execution pursuit coalition to capture an evader. The second inequality constraint indicates that each evader is assigned at most one pursuit coalition, and the third one restricts that each pursuer can occur at most one pursuit coalition when the game runs. It can be verified that \autoref{matrix} can find all execution pursuit coalitions which contain a specific pursuer. Thus, this 0-1 integer programming is solvable. \qed

\begin{algm}\label{matrix}\rm
(Matrix for Inequality Constraint).\\
1: \textbf{Input:} $N_p,N_e$, and $A_3\gets$null matrix.\\
2: \textbf{if} $N_p\ge2$ \textbf{then}\\
3: \textbf{for} $i\in\{1,2,...,N_p\}$ \textbf{do}\\
4: \hspace{0.2 in}$\rm temp\gets$null matix\\
5: \hspace{0.2 in}\textbf{for} $j\in\{1,2,...,N_p(N_p-1)/2\}$ \textbf{do}\\
6: \hspace{0.4 in}$\rm tag\gets0$\\
7: \hspace{0.4 in}\textbf{for} $k\in\{1,2,...,i-1\}(i\ge2)$ \textbf{do}\\
8: \hspace{0.6 in}\textbf{if} $j==i-k+(k-1)(N_p-k/2)$ \textbf{then}\\
9: \hspace{0.77 in}$\rm tag\gets1$,\textbf{break} \textbf{ end if end for}\\
10:\hspace{0.37 in}\textbf{if} $j\ge(i-1)(N_p-i/2)+1$ \&\&\\
11:\hspace{0.51 in}$j\leq i(N_p-(i+1)/2)$ \textbf{then}\\
12:\hspace{0.57 in}$\rm tag\gets1$ \textbf{end if}\\
13:\hspace{0.37 in}$\rm temp\gets[temp,tag]$ \textbf{end for}\\
14:\hspace{0.17 in}$A_3\gets[A_3;\rm temp]$ \textbf{end for}\\
15: $A_3\gets A_3\otimes {\rm{ones}}(1,N_e)$ \textbf{end if}\\\
16: $A_3\gets [I_{N_p}\otimes {\rm ones}(1,N_e),A_3]$\\
17:\textbf{Output:} $A_3$.
\end{algm}

\begin{rek}\label{rek43}\rm
It can be observed that $q$ is unique, while the maximum matching $\bm{z}^*$ may have multiple solutions. Moreover, the original matching is a $(2^{N_p}-1)\times N_e$ constrained bipartite matching problem, which is a NP-problem. However, the maximum matching given by \autoref{maximum} is a P-problem with $N_v$ variables and $N_v+N_e+N_p$ inequality constraints. The first inequality constraint in (\ref{maximizeequal11}) will be simplified a lot when $\bm{r}$ contains many zero elements.
\end{rek}

\begin{defi}[Maximum Matching Pairs]\label{pairs}\rm
For a maximum matching $\bm{z}^*=[\bm{s}^{1*},\bm{s}^{2*}]^\mathsf{T}$, define the following sets of matching pairs:\begin{equation}\begin{aligned}
 M^1(\bm{z}^*)&=\big\{(i,j)|s_i^{1*}(j)=1,1\leq i\leq N_p,1\leq j\leq N_e\big\}\\
 M^2(\bm{z}^*)&=\big\{(i1,i2,j)|s_{i1,i2}^{2*}(j)=1,1\leq i1<i2\leq N_p,\ \ \\
&\qquad\qquad\qquad\qquad\qquad\qquad\qquad\ 1\leq j\leq N_e\big\}.
\end{aligned}\end{equation}
\end{defi}

\begin{figure}\centering
\graphicspath{{Figures/}}
\subfigure{
\label{fig:sub-1}
\includegraphics[width=40mm,height=33mm]{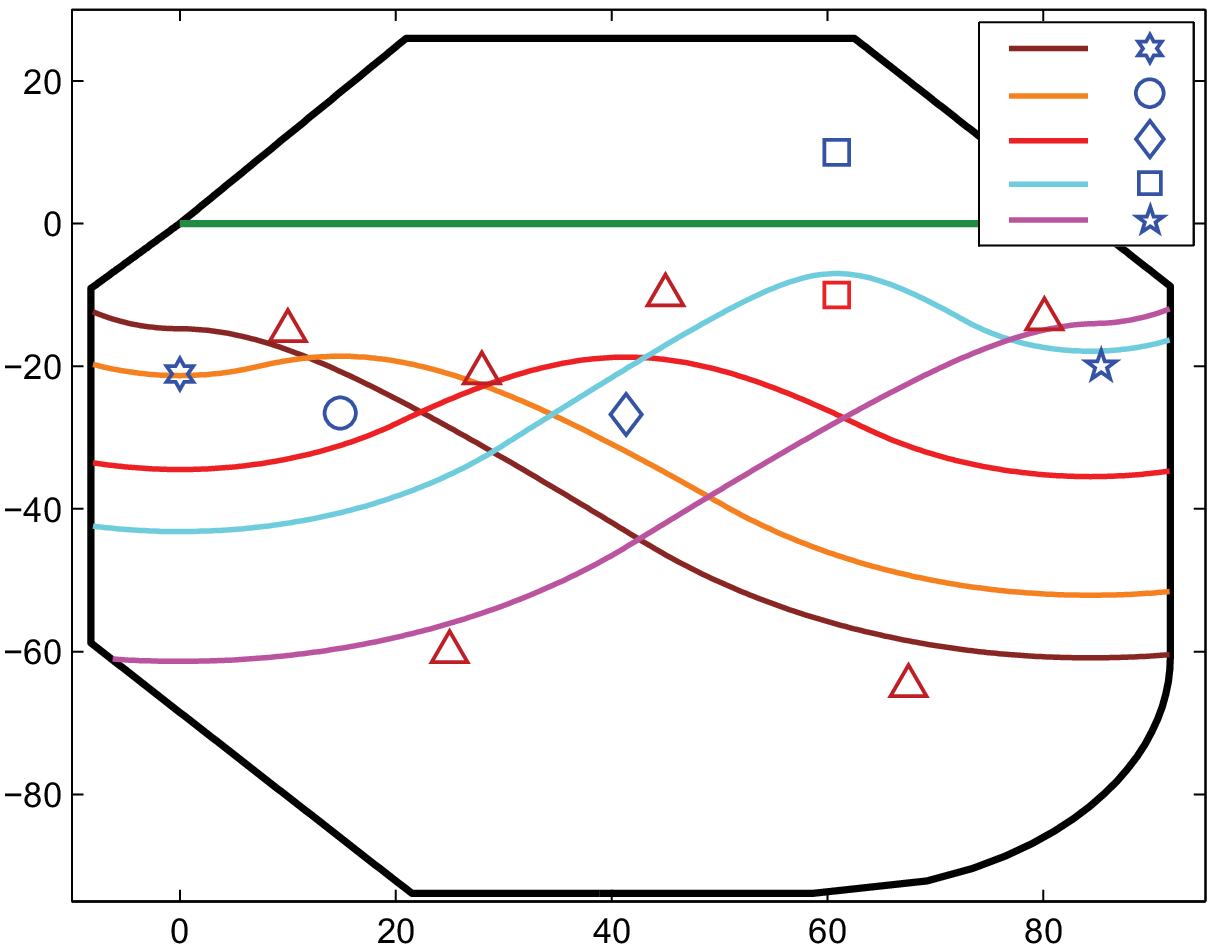}}
\put(-75,9){$\scriptstyle{\Omega_{\rm play}}$}
\put(-75,83){$\scriptstyle{\Omega_{\rm tar}}$}
\put(-75,65){$\scriptstyle{\mathcal{T}}$}
\put(-60,-8){$\scriptstyle{(a)}$}
\subfigure{
\label{fig:sub-2}
\includegraphics[width=41mm,height=33mm]{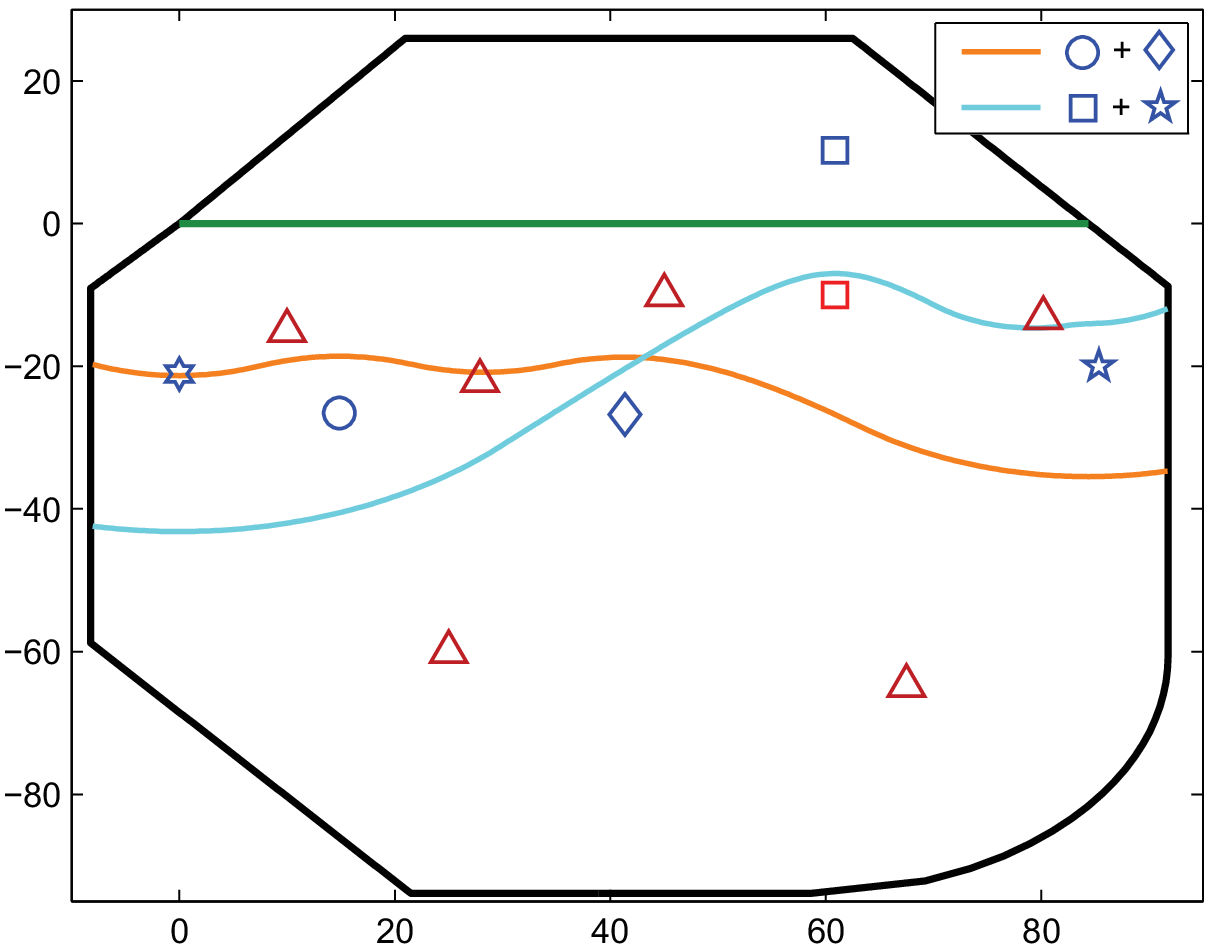}}
\put(-76,9){$\scriptstyle{\Omega_{\rm play}}$}
\put(-76,83){$\scriptstyle{\Omega_{\rm tar}}$}
\put(-76,65){$\scriptstyle{\mathcal{T}}$}
\put(-60,-8){$\scriptstyle{(b)}$}\\
\subfigure{
\label{fig:sub-3}
\includegraphics[width=40mm,height=33mm]{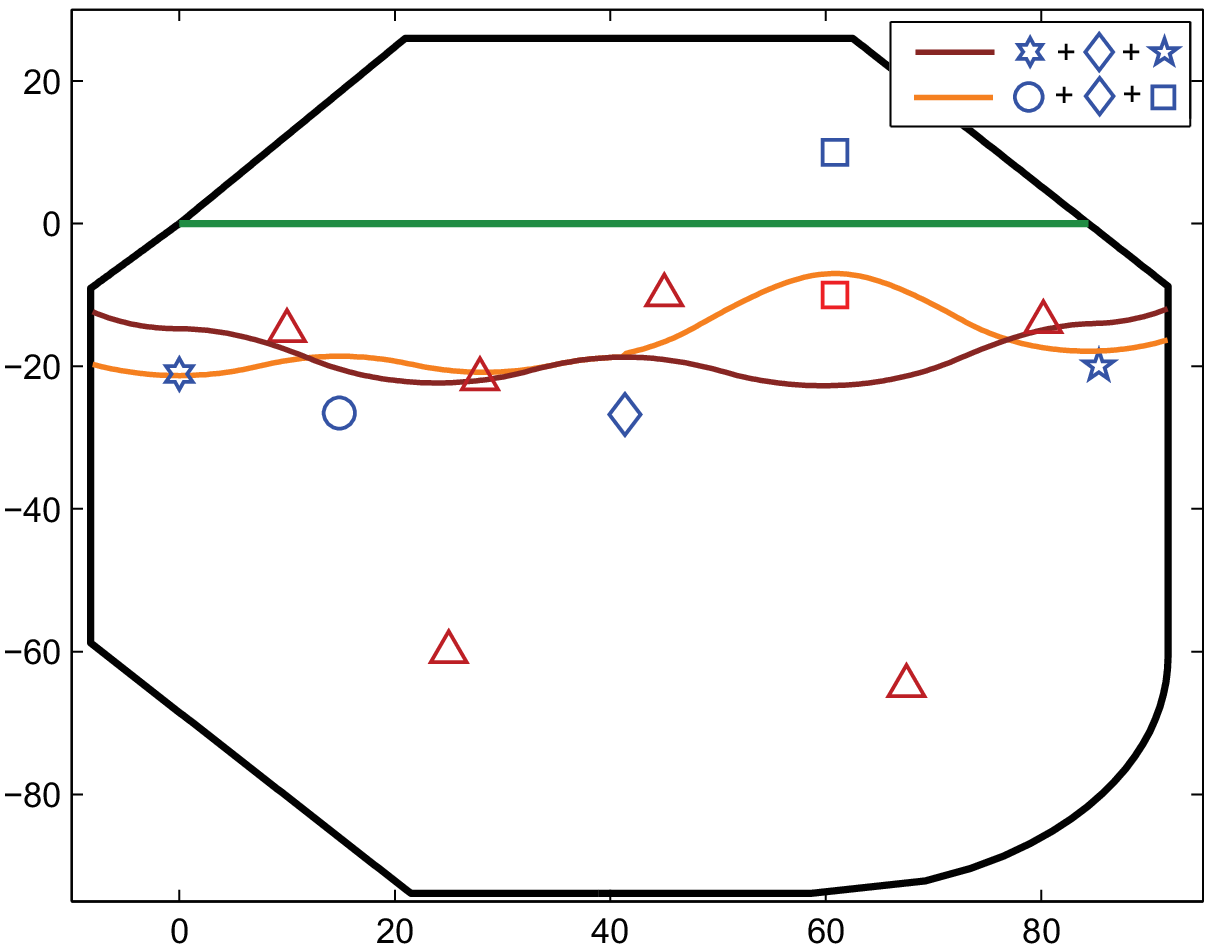}}
\put(-75,9){$\scriptstyle{\Omega_{\rm play}}$}
\put(-75,83){$\scriptstyle{\Omega_{\rm tar}}$}
\put(-75,65){$\scriptstyle{\mathcal{T}}$}
\put(-60,-8){$\scriptstyle{(c)}$}
\subfigure{
\label{fig:sub-4}
\includegraphics[width=41mm,height=33mm]{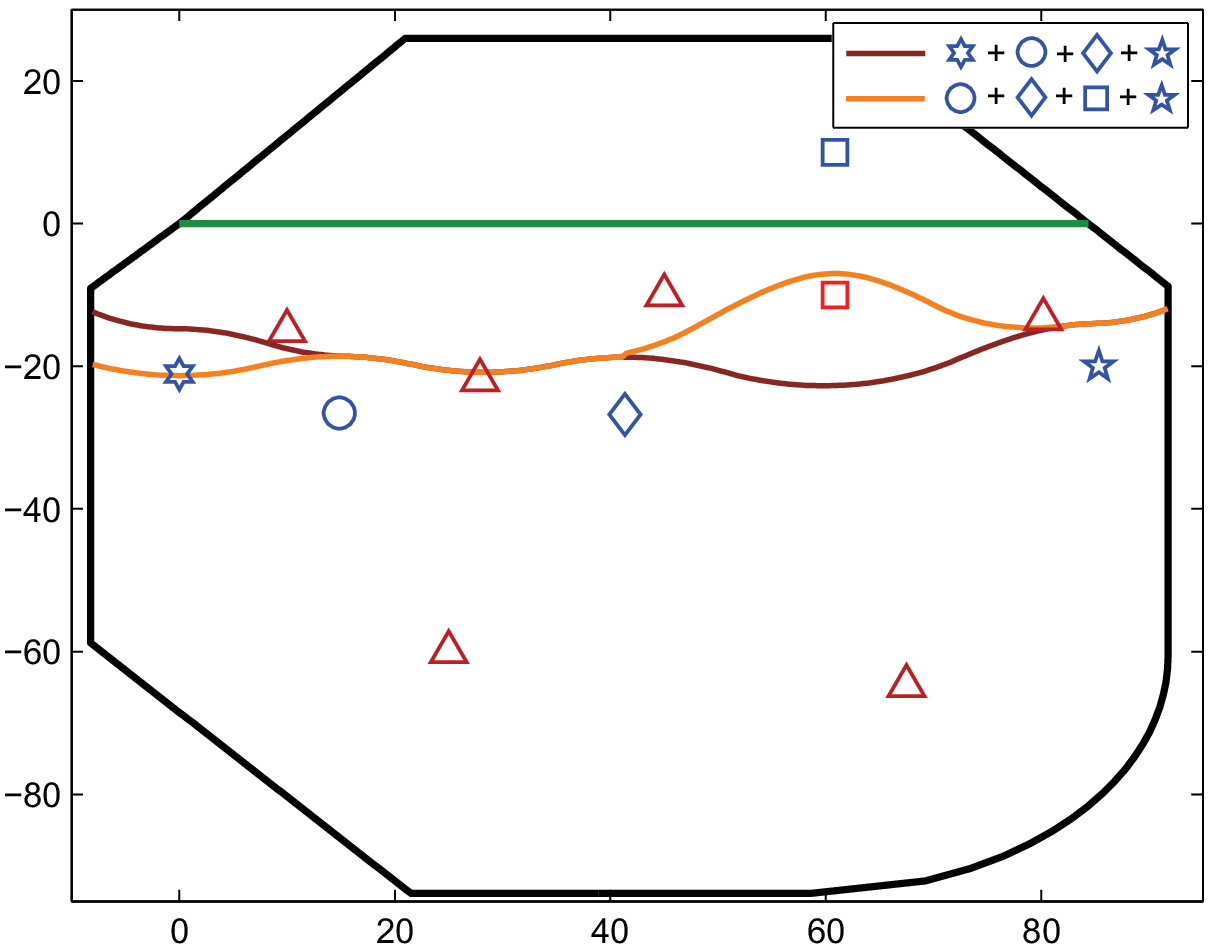}}
\put(-76,9){$\scriptstyle{\Omega_{\rm play}}$}
\put(-76,83){$\scriptstyle{\Omega_{\rm tar}}$}
\put(-76,65){$\scriptstyle{\mathcal{T}}$}
\put(-60,-8){$\scriptstyle{(d)}$}\\
\begin{flushleft}\ \ 
\subfigure{
\label{fig:sub-5}
\includegraphics[width=40mm,height=33mm]{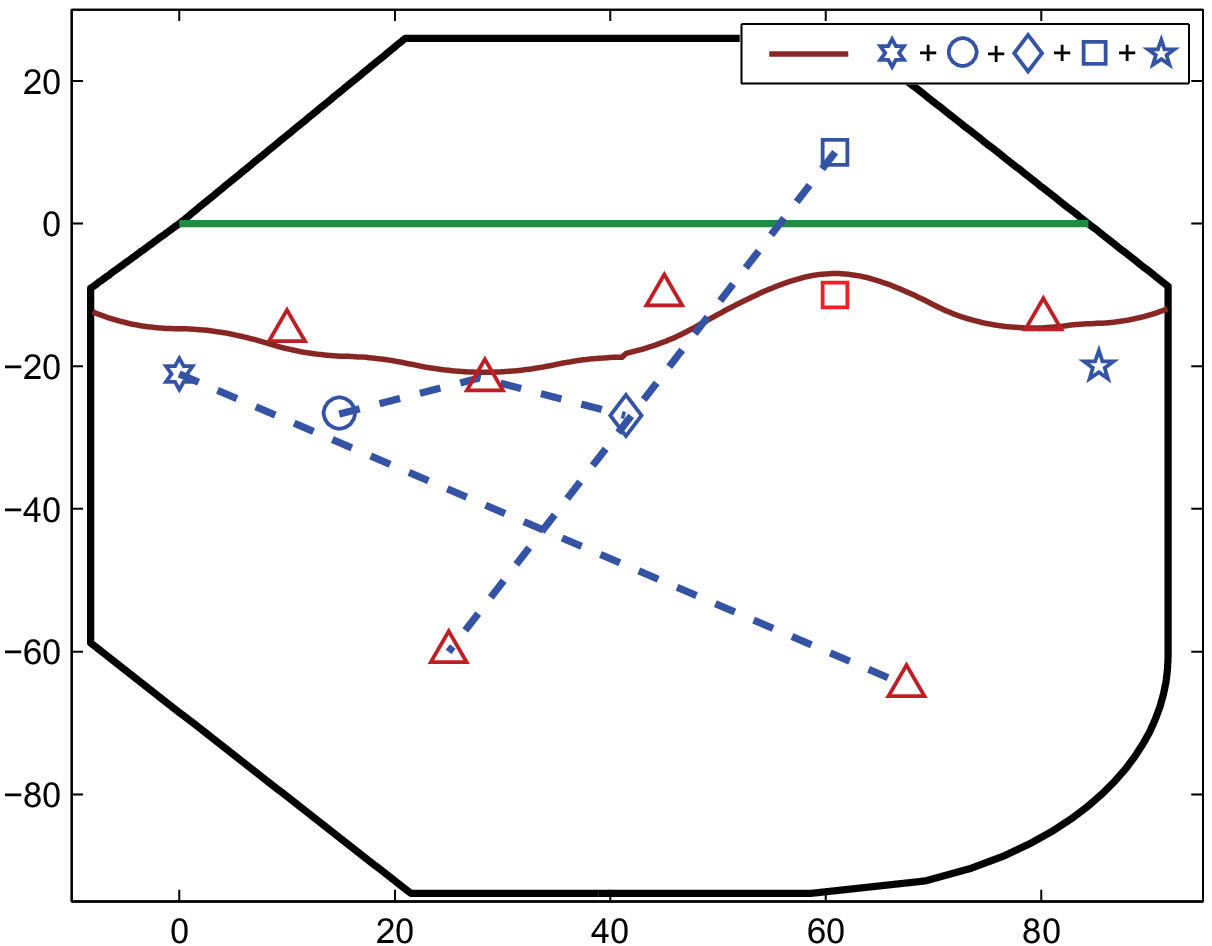}}
\put(-75,9){$\scriptstyle{\Omega_{\rm play}}$}
\put(-75,83){$\scriptstyle{\Omega_{\rm tar}}$}
\put(-75,65){$\scriptstyle{\mathcal{T}}$}
\put(-60,-8){$\scriptstyle{(e)}$}
\subfigure{
\label{fig:sub-6}
\includegraphics[width=30mm,height=33mm]{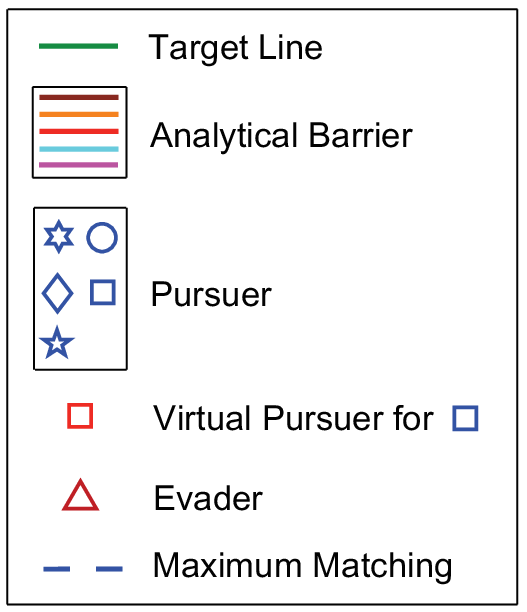}}
\end{flushleft}
\caption{5 pursuers versus 6 evaders. $(a)-(e)$ Computation of the barrier and winning regions for all pursuit coalitions when facing one evader: $(a)$ all one-pursuer pursuit coalitions; $(b)$ two-pursuer pursuit coalitions (only two are depicted); $(c)$ three-pursuer pursuit coalitions (only two are depicted); $(d)$ four-pursuer pursuit coalitions (only two are depicted); $(5)$ five-pursuer pursuit coalition. Fig. \ref{fig:sub-5} also shows the generated maximum matching: two one-to-one matching pairs and one two-to-one matching pair. Thus, the pursuit team currently can guarantee to capture at most three evaders in $\Omega_{\rm play}$.}
\label{figureexample1}
\end{figure}

Note that $M^1(\bm{z}^*)$ represents all one-to-one matching pairs in $\bm{z}^*$, and $M^2(\bm{z}^*)$ represents all two-to-one matching pairs.

\section{Simulation Results}\label{simulation}
In this section, simulation results are presented to illustrate the previous theoretical developments. Assume that $v_P=1{\rm m/s},v_E=0.7{\rm m/s}$, namely, $\alpha=0.7$, $N_p=5$ and $N_e=6$.

Fig. \ref{figureexample1} shows the barrier and winning regions for all pursuit coalitions, namely, their capturable and uncapturable regions, where Fig. \ref{fig:sub-1} and Fig. \ref{fig:sub-5} refer to all one-pursuer pursuit coalitions and five-pursuer pursuit coalition, respectively. For clarity, only two barriers and winning regions are depicted for two-pursuer, three-pursuer, and four-pursuer pursuit coalitions in Fig. \ref{fig:sub-2}, Fig. \ref{fig:sub-3}, and Fig. \ref{fig:sub-4}, respectively.

Fig. \ref{fig:sub-5} also shows the maximum matching from the known barrier and winning regions, including two one-to-one matching pairs and one two-to-one matching pair. Thus, the maximum matching is of size 3. This guarantees that if each pursuer occurring in the maximum matching plays optimally against the evader matched by the maximum matching, this matched evader will be captured before reaching $\mathcal{T}$.

\section{Conclusions and Future Work}\label{conclusion}
\subsection{Conclusions}
This paper considered a multiplayer reach-avoid game in a general convex domain. The key achievement is providing an analytical description of the winning regions when each possible pursuit coalition competes with one evader. Furthermore, an interception scheme involving pursuer-evader matching is generated for the pursuit team such that the most evaders can be captured before entering the target region. The constructed barrier, splitting the winning regions, shows that at most two pursuers are needed to intercept one evader if the capture is possible, greatly simplifying the matching search. 

The winning regions from the whole pursuit team and one evader case also can help the evasion team to determine which evaders can enter the target region or escape from the play region definitely, no matter what strategy the pursuit team uses. Then these evaders will get more attention from the evasion team and be chosen as key mission performers. More generally, this result can provide guarantees on goal satisfaction and safety of optimal system trajectories for safety-critical systems, where a group of vehicles aim to reach their destinations in the presence of dynamic obstacles. 

The results in this paper are almost analytical and applicable for real-time updates. More importantly, all possible pursuit coalitions are considered and our results are optimal in the sense of the most evaders intercepted.

\subsection{Future Work}
All players considered in this paper move with simple motion and are able to turn instantaneously. In future work, more practical and complex dynamic models will be considered, for example, those of the Isaacs-Dubins car. Extending the results obtained in general convex domains to nonconvex domains is another worth-pursuing direction, and extracting out some conservative results is straightforward, such as, approximating the nonconvex domain with an appropriate convex domain. In view of the limited communication and computing power of a single player, another interesting and promising possibility for future extension is to consider distributed multiplayer reach-avoid games, which are the focus of our following research.


\bibliographystyle{IEEEtran}
\bibliography{Bibliography/BIB_1x-TIE-2xxx}

\end{document}